\documentclass[12pt]{article} 
\usepackage{latexsym, 
amscd, 
graphicx, epsfig, amssymb}

\newtheorem{thm}{Theorem}[section]
\newtheorem{defn}[thm]{Definition}
\newtheorem{prop}[thm]{Proposition}
\newtheorem{cor}[thm]{Corollary}
\newtheorem{lemma}[thm]{Lemma}
\newtheorem{rema}[thm]{Remark}

\newcommand{\halmos}{\rule{1ex}{1.4ex}}

\newcommand{\text}[1]{\mbox{\rm #1}}

\newcommand{\nn}{\nonumber \\}

 \newcommand{\res}{\mbox{\rm Res}}

 \newcommand{\pf}{{\it Proof.}\hspace{2ex}}
 \newcommand{\epf}{\hspace*{\fill}\mbox{$\halmos$}}
 \newcommand{\epfv}{\hspace*{\fill}\mbox{$\halmos$}\vspace{1em}}

\newcommand{\wt}{\mbox{\rm wt}\ }
\newcommand{\swt}{\mbox{\rm {\scriptsize wt}}\ }

\newcommand{\Y}{\mathcal{Y}}
\newcommand{\C}{\mathbb{C}}
\newcommand{\Z}{\mathbb{Z}}
\newcommand{\R}{\mathbb{R}}
\newcommand{\Q}{\mathbb{Q}}
\newcommand{\N}{\mathbb{N}}

\newcommand{\one}{\mathbf{1}}

\makeatletter
\newlength{\@pxlwd} \newlength{\@rulewd} \newlength{\@pxlht}
\catcode`.=\active \catcode`B=\active \catcode`:=\active \catcode`|=\active
\def\sprite#1(#2,#3)[#4,#5]{
   \edef\@sprbox{\expandafter\@cdr\string#1\@nil @box}
   \expandafter\newsavebox\csname\@sprbox\endcsname
   \edef#1{\expandafter\usebox\csname\@sprbox\endcsname}
   \expandafter\setbox\csname\@sprbox\endcsname =\hbox\bgroup
   \vbox\bgroup
      \catcode`.=\active\catcode`B=\active\catcode`:=\active\catcode`|=\active
      \@pxlwd=#4 \divide\@pxlwd by #3 \@rulewd=\@pxlwd
      \@pxlht=#5 \divide\@pxlht by #2
      \def .{\hskip \@pxlwd \ignorespaces}
      \def B{\@ifnextchar B{\advance\@rulewd by \@pxlwd}{\vrule
         height \@pxlht width \@rulewd depth 0 pt \@rulewd=\@pxlwd}}
      \def :{\hbox\bgroup\vrule height \@pxlht width 0pt depth 0pt\ignorespaces}
      \def |{\vrule height \@pxlht width 0pt depth 0pt\egroup
         \prevdepth= -1000 pt}
   }
\def\endsprite{\egroup\egroup}
\catcode`.=12 \catcode`B=11 \catcode`:=12 \catcode`|=12\relax
\makeatother

\def\hboxtr{\FormOfHboxtr} 
\sprite{\FormOfHboxtr}(25,25)[0.5 em, 1.2 ex] 

:BBBBBBBBBBBBBBBBBBBBBBBBB |
:BB......................B |
:B.B.....................B |
:B..B....................B |
:B...B...................B |
:B....B..................B |
:B.....B.................B |
:B......B................B |
:B.......B...............B |
:B........B..............B |
:B.........B.............B |
:B..........B............B |
:B...........B...........B |
:B............B..........B |
:B.............B.........B |
:B..............B........B |
:B...............B.......B |
:B................B......B |
:B.................B.....B |
:B..................B....B |
:B...................B...B |
:B....................B..B |
:B.....................B.B |
:B......................BB |
:BBBBBBBBBBBBBBBBBBBBBBBBB |

\endsprite

\title{ {\bf Cofiniteness conditions,
projective covers and the logarithmic tensor product theory} }
\date{}
\author{Yi-Zhi Huang}

\begin{document}

\bibliographystyle{alpha}
\maketitle

\begin{abstract} 
We construct projective covers of irreducible $V$-modules 
in the category of
grading-restricted generalized $V$-modules
when $V$ is a vertex operator algebra satisfying 
the following conditions: 1. $V$ is $C_{1}$-cofinite in the sense of Li.
2. There exists a positive integer $N$ such that 
the differences between the real parts of the lowest conformal
weights of irreducible $V$-modules are bounded by $N$
and such that the associative algebra
$A_{N}(V)$ is finite dimensional. This result 
shows that the category of
grading-restricted generalized $V$-modules is a finite abelian 
category over $\C$.
Using the 
existence of projective covers, we prove that if such a vertex operator 
algebra $V$ satisfies in addition Condition~3,
that irreducible $V$-modules are $\R$-graded and $C_{1}$-cofinite
in the sense of the author,
then the category of grading-restricted
generalized $V$-modules is closed under operations 
$\hboxtr_{P(z)}$ for $z\in \C^{\times}$. We also prove that 
other conditions for applying the logarithmic tensor product 
theory developed by Lepowsky, Zhang and the author hold. 
Consequently, for such $V$, 
this category has a natural structure 
of braided tensor category. In particular, when $V$ is of positive 
energy and $C_{2}$-cofinite, Conditions 1--3
are satisfied and thus all the conclusions hold.
\end{abstract}

\renewcommand{\theequation}{\thesection.\arabic{equation}}
\renewcommand{\thethm}{\thesection.\arabic{thm}}
\setcounter{equation}{0}
\setcounter{thm}{0}
\setcounter{section}{-1}

\section{Introduction}

In the present paper, we construct projective covers of irreducible
$V$-modules in the category of grading-restricted generalized
$V$-modules and prove that the logarithmic tensor product theory
developed in \cite{HLZ1} and \cite{HLZ2} can be applied to this
category when $V$ satisfies suitable natural 
cofiniteness and other conditions
(see below). Consequently, for such $V$, 
this category is a finite abelian category over $\C$ 
and has a natural structure 
of braided tensor category. We refer the reader to \cite{EO}, 
\cite{Fu} and
\cite{HLZ2} for detailed discussions on the importance 
and applications of projective covers and logarithmic tensor 
products.

For a vertex operator algebra $V$ satisfying certain finite reductivity
conditions, Lepowsky and the author have developed a tensor product theory
for $V$-modules in \cite{tensor0}--\cite{tensor5}, \cite{H1} and
\cite{H6}.  Consider a simple vertex operator algebra $V$ satisfies the
following slightly stronger conditions: (i) $V$ is of positive energy
(that is, $V_{(n)}=0$ for $n<0$ and $V_{(0)}=\C\one$) and $V'$ is equivalent 
to $V$ as $V$-modules. (ii) Every
$\N$-gradable weak $V$-module is completely reducible. (iii) $V$ is
$C_{2}$-cofinite.  Then the author further proved in \cite{H11} (see
also \cite{H8}, \cite{H9}; cf. \cite{Le})
that the category of $V$-modules for such $V$ has a natural structure of
modular tensor category.  This result reduces a large part of the
representation theory of such a vertex operator algebra to the study of
the corresponding
modular tensor category and allows us to employ the powerful
homological-algebraic methods. 

The representation theory of a vertex operator algebra satisfying the
three conditions above corresponds to a chiral conformal field 
theory\footnote{In this introduction, though we often 
mention conformal
field theories, we shall not discuss the precise mathematical formulation
of conformal field theory and the problem of mathematically constructing 
conformal field theories. Instead, we use the term conformal field theories 
to mean certain conformal-field-theoretic structures and results, 
such as operator product expansions, modular invariance, fusion rules,
Verlinde formula and so on. See \cite{H2}--\cite{H10} and 
\cite{HK1}--\cite{HK2} for the relationship
between the representation theory of vertex operator algebras and conformal 
field theories and the results on the mathematical formulation and 
construction of conformal field theories in terms of representations of 
vertex operator algebras.}.
Such a chiral conformal field theory has
all the properties of the chiral part of 
a rational conformal field theory. In fact, in view of
the results in \cite{H10} and \cite{H11} (see also \cite{H8}, \cite{H9};
cf. \cite{Le}), one might even want to 
define a rational conformal field theory to be a conformal field
theory whose chiral algebra is a vertex operator algebra satisfying the
three conditions above (see, for example, \cite{Fu}). 

In the study of many problems in mathematics and physics, 
for example, problems in the studies of 
mirror symmetry, string theory, disorder systems
and percolation, it is necessary to study irrational conformal 
field theories. If we use the definition above as the definition 
of rational conformal field theory, then to study irrational 
conformal field theories means that we have to study the 
representation theory of vertex operator algebras for which
at least one of the three conditions is not satisfied. 

In the present paper, 
we study the representation theory of vertex
operator algebras satisfying certain  conditions
weaker than Conditions (i)--(iii) above. In particular, we shall not 
assume the complete reducibility of $\N$-gradable $V$-modules 
or even $V$-modules. 
Since the complete reducibility is not assumed, one will not be able
to generalize the author's proof in \cite{H6} to show that in this case
the analytic extensions of products of intertwining operators still do
not have logarithmic terms.  Thus in this case, the corresponding
conformal field theories in physics must be logarithmic ones
studied first by Gurarie in \cite{G}. The triplet $\mathcal{W}$-algebras
of central charge $1-6\frac{(p-1)^{2}}{p}$,
introduced first by Kausch \cite{K1} and studied extensively
by Flohr \cite{F1} \cite{F2}, 
Kausch \cite{K2}, Gaberdiel-Kausch \cite{GK1} \cite{GK2}, 
Fuchs-Hwang-Semikhatov-Tipunin \cite{FHST}, Abe \cite{A},
Feigin-Ga{\u\i}nutdinov-Semikhatov-Tipunin \cite{FGST1} \cite{FGST2}
\cite{FGST3}, Carqueville-Flohr \cite{CF}, Flohr-Gaberdiel 
\cite{FG}, Fuchs \cite{Fu} and Adamovi\'{c}-Milas \cite{AM1} \cite{AM2},
are examples of vertex operator algebras satisfying the positive energy 
condition and 
the $C_{2}$-cofiniteness condition, but not 
Condition (ii) above. For the proof of 
the $C_{2}$-cofiniteness condition, see \cite{A} for the simplest
$p=2$ case and \cite{CF} and \cite{AM2} for the general case. 
For the proof that Condition (ii) is not satisfied by these
vertex operator algebras, see \cite{A} for the simplest  $p=2$ case
and \cite{FHST} and \cite{AM2} for the general case. 
A family of $N=1$ triplet vertex operator superalgebras 
has been constructed and studied recently by Adamovi\'{c}
and Milas in \cite{AM3}. Among many results obtained in \cite{AM3} 
are the $C_{2}$-cofiniteness of these vertex operator superalgebras 
and a proof that Condition (ii) is not satisfied by 
them.

In \cite{HLZ1} and \cite{HLZ2}, Lepowsky, Zhang and the author
generalized the tensor product theory of Lepowsky and the author
\cite{tensor0}--\cite{tensor5} \cite{H1} \cite{H6} to a logarithmic
tensor product theory for suitable categories of generalized modules for
M\"{o}bius or conformal vertex algebras satisfying suitable conditions.
In this theory, generalized modules in these categories are not required
to be completely reducible, not even required to be completely reducible
for the operator $L(0)$ (see also \cite{M}). 
The general theory in \cite{HLZ1} and
\cite{HLZ2} is quite flexible since it can be applied to any category of
generalized modules such that the assumptions in \cite{HLZ1} and
\cite{HLZ2} hold.  One assumption is that the
category should be closed under the $P(z)$-tensor product $\boxtimes_{P(z)}$
for some $z\in \C^{\times}$. Since the category is also assumed to be 
closed under the operation of taking contragredient, this assumption
is equivalent to the assumption that the category is closed under
an operation $\hboxtr_{P(z)}$ (see \cite{HLZ1} and \cite{HLZ2}). There are also
some other assumptions for the category to be a braided
tensor category. 

This logarithmic tensor product theory can be applied to 
a range of different examples. 
The original tensor product theory 
developed in \cite{tensor0}--\cite{tensor5}, \cite{H1} and \cite{H2}
becomes a special case. We also expect that this logarithmic 
tensor product theory will play an important role
in the study of unitary conformal field theories which 
do not have logarithmic fields but are 
not necessarily rational. 
For a vertex operator algebra associated to modules for an affine Lie 
algebra of a non-positive integral level, Zhang \cite{Zha1} \cite{Zha2}
proved that the category $\mathcal{O}_{\kappa}$ 
is closed under the operation $\hboxtr_{P(z)}$ 
by reinterpreting,  in the framework of 
\cite{HLZ1} and \cite{HLZ2}, the result proved by Kazhdan and Lusztig 
\cite{KL1}--\cite{KL5} that the category $\mathcal{O}_{\kappa}$ 
is closed under 
their tensor product bifunctor. 
It is also easy to see that objects in the category $\mathcal{O}_{\kappa}$
satisfy the $C_{1}$-cofiniteness condition in the 
sense of \cite{H6} (see \cite{Zha2})
and the category $\mathcal{O}_{\kappa}$ 
satisfies the other conditions for applying the logarithmic tensor 
product theory in \cite{HLZ1} and \cite{HLZ2}. 
As a consequence, 
we obtain another construction of the Kazhdan-Lusztig braided 
tensor category structures.

In the present paper, we consider the following conditions 
for a vertex operator algebra
$V$: 1. $V$ is $C_{1}$-cofinite in the sense of Li \cite{L}.
2. There exists a positive integer $N$ such that 
the differences of the real parts of the lowest 
conformal weights of irreducible $V$-modules are less than or equal
to $N$
and such that the associative algebra
$A_{N}(V)$ introduced by  
Dong-Li-Mason \cite{DLM1} (a generalization of the associative 
algebra $A_{0}(V)$ introduced by Zhu \cite{Zhu2}) is finite dimensional.
3. Irreducible $V$-modules are $\R$-graded and
are  $C_{1}$-cofinite in the 
sense of \cite{H6}. Note that the first part of 
Condition 2 is always satisfied when $V$ has only finitely many 
irreducible $V$-modules and also note that if $V$ is of positive energy and 
$C_{2}$-cofinite, $V$ satisfies all the three conditions 1--3
(see Proposition \ref{c-2+p-e=>cond-1-3}).

When $V$ satisfies Conditions 1 and 2, we prove that 
a generalized $V$-module is of finite-length if and only 
if it is grading restricted and if and only if it is 
quasi-finite dimensional.
Grading-restricted generalized $V$-modules
were first studied by Milas in \cite{M}
and were called logarithmic 
$V$-modules.  When $V$ satisfies these two conditions, 
we prove, among many 
basic properties of generalized $V$-modules, that  
any irreducible $V$-module $W$ has a projective cover
in the category of
grading-restricted generalized $V$-modules.
The existence of projective covers 
of irreducible $V$-modules
is a basic assumption in 
\cite{Fu} and, since Condition 2 implies that there are only finitely many 
irreducible $V$-modules, this existence says that the category of 
grading-restricted generalized $V$-modules is a finite abelian 
category over $\C$.

Using this 
existence of projective covers, we prove that if
$V$ satisfies Conditions 1--3,
the category of grading-restricted
generalized $V$-modules is closed under the operation 
$\hboxtr_{P(z)}$ for $z\in \C^{\times}$ and thus is 
closed under the operation $\boxtimes_{P(z)}$. 
We also prove that 
other conditions for applying the logarithmic tensor product 
theory in \cite{HLZ1} and \cite{HLZ2}
developed by Lepowsky, Zhang and the author hold. 
Consequently,  this category has a natural structure of a 
braided tensor category.

Note that if, in addition, $V$ is simple and the braided tensor category
of grading-restricted generalized $V$-modules is rigid, then the results
of the present paper show that this category is a ``finite tensor
category" in the sense of Etingof and Ostrik \cite{EO} (cf. \cite{Fu}). 
We shall discuss the rigidity in a future paper since
it needs the generalization of the author's
result \cite{H7} on the modular invariance of the $q$-traces of 
intertwining operators in the finitely reductive case to
a result on the modular invariance of the $q$-(pseudo-)traces of 
logarithmic intertwining operators in the nonreductive case.

Since a positive energy
$C_{2}$-cofinite vertex operator algebra satisfies
Conditions 1--3, these conclusions hold for such a 
vertex operator algebra\footnote{In a preprint arXiv:math/0309350,
incorrect ``tensor product" operations in 
the category of finite length generalized
modules for a vertex operator algebra satisfying the 
$C_{2}$-cofiniteness condition were introduced. The counterexamples 
given in the paper \cite{HLLZ} show
that the ``tensor product"
operations in that 
preprint do not have the basic properties a tensor product 
operation must have and are different from the tensor product 
operations defined in \cite{HLZ1} and \cite{HLZ2} and 
proved to exist in the present paper (see the paper \cite{HLLZ}
for details).}.

For many results in this paper, we prove them under weaker assumptions
so that these results might also be useful for other purposes.  Some of
the proofs can be simplified greatly if $V$ satisfies stronger
conditions such as the $C_{2}$-cofiniteness condition but we shall 
not discuss these simplifications.  Almost all the
results in the present paper hold also when the vertex operator
algebra is a grading-restricted M\"{o}bius vertex algebra.

The present paper is organized as follows: In Section 1, we 
given basic definitions and properties for generalized 
modules. We study cofiniteness conditions for vertex operator 
algebras and their modules in Section 2. We also discuss  
associative algebras introduced by Zhu and Dong-Li-Mason
in this section. In Section 3,  we prove that if $V$ satisfies 
Conditions 1 and 2 above, then 
any irreducible $V$-module $W$ has a projective cover
in the category of 
grading-restricted generalized modules. In particular, 
we see that this category is a finite abelian category over $\C$.
In Section 4, we prove that if $V$ satisfies Conditions 1--3,
the category of grading-restricted
generalized $V$-modules is closed under the operation $\hboxtr_{P(z)}$
for any $z\in \mathbb{Z}^{\times}$ and thus is closed under the 
$P(z)$-tensor product $\boxtimes_{P(z)}$. The 
other assumptions needed in the logarithmic tensor product 
theory are also shown to hold in this section. 
Combining with the results of \cite{HLZ1} and \cite{HLZ2}, 
we obtain the conclusion that this category is a braided tensor 
category.

\paragraph{Acknowledgment} 
The author gratefully acknowledges the partial support
from NSF grant DMS-0401302. The author is grateful to 
Antun Milas and Jim Lepowsky for comments. 

\renewcommand{\theequation}{\thesection.\arabic{equation}}
\renewcommand{\thethm}{\thesection.\arabic{thm}}
\setcounter{equation}{0}
\setcounter{thm}{0}

\section{Definitions and basic properties}

In this paper, we shall assume that the reader is familiar with the basic 
notions and results in the theory of vertex operator algebras.
In particular, we assume that the reader is 
familiar with weak modules,
$\N$-gradable weak modules, contragredient modules and 
related results. Our terminology and conventions follow
those in \cite{FLM}, \cite{FHL} and \cite{LL}. 
We shall use $\Z$, $\Z_{+}$, $\N$, $\Q$, $\R$, 
$\C$ and $\C^{\times}$ to denote the (sets of) integers, positive integers,
nonnegative integers, rational numbers, real numbers, complex numbers
and nonzero complex numbers, respectively. For $n\in \C$, we use
$\Re{(n)}$ and $\Im{(n)}$ to denote the real and imaginary parts of 
$n$.

We fix a vertex operator algebra $(V, Y, \one, 
\omega)$ in this paper. We first recall the definitions of 
generalized $V$-module and related notions in \cite{HLZ1}
and \cite{HLZ2} (see also \cite{M}):

\begin{defn}
{\rm A {\it generalized $V$-module} is a $\C$-graded 
vector space $W=\coprod_{n\in \C}W_{[n]}$ equipped with 
a linear map 
\begin{eqnarray*}
Y_{W}: V\otimes W&\to &W((x))\\
v&\mapsto & Y_{W}(v, x)
\end{eqnarray*}
satisfying all the axioms for $V$-modules except
that we do not require $W$ satisfying the two grading-restriction
conditions and that 
the $L(0)$-grading property is replaced by the following weaker version,
still called the {\it $L(0)$-grading property}:
For $n\in \C$, the homogeneous
subspaces $W_{[n]}$ is the generalized eigenspaces of $L(0)$ with 
eigenvalues $n$, that is, for $n\in \C$, 
there exists $K\in \Z_{+}$ such that 
$(L(0)-n)^{K}w=0$.
{\it Homomorphisms} (or {\it module maps})
and {\it isomorphisms} (or {\it equivalence})
between generalized $V$-modules, {\it generalized $V$-submodules}
and {\it quotient generalized $V$-modules} 
are defined in the obvious way.}
\end{defn}

The generalized modules we are mostly interested in the present paper are 
given in the following definition:

\begin{defn}
{\rm A generalized $V$-module $W$ is {\it irreducible} 
if there is no generalized $V$-submodule of $W$ which is neither
$0$ nor $W$ itself.  A generalized $V$-module is {\it lower truncated}
if $W_{[n]}=0$ when $\Re{(n)}$ is sufficiently negative.
For a lower-truncated generalized $V$-module $W$, 
if there exists $n_{0}\in \C$ such that 
$W_{[n_{0}]}\ne 0$ but $W_{[n]}=0$ when $\Re{(n)}<\Re{(n_{0})}$ 
or $\Re{(n)}=\Re{(n_{0})}$ but $\Im{(n)}\ne \Im{(n_{0})}$,
then we say that $W$ {\it has a lowest conformal weight},
or for simplicity, $W$ {\it has a lowest weight}. In this case, 
$n_{0}$ is 
called the {\it lowest conformal weight} or {\it lowest 
weight} of $W$, the homogeneous subspace 
$W_{[n_{0}]}$ of $W$ is called the {\it lowest weight space}
or {\it lowest weight space}
of $W$ and elements of $W_{[n_{0}]}$ are called {\it lowest 
conformal weight 
vectors} or {\it lowest weight vectors} of $W$.
A generalized $V$-module is {\it grading restricted} if $W$
is lower truncated 
and $\dim W_{[n]}<\infty$
for $n\in \C$. 
A {\it quasi-finite-dimensional generalized $V$-module}
is a generalized $V$-module such that for any real number 
$R$, $\dim \coprod_{\Re{(n)}\le R}W_{[n]}<\infty$.
A generalized $V$-module $W$ is
an {\it (ordinary) $V$-module} if $W$ is grading restricted and 
$W_{[n]}=W_{(n)}$ for $n\in \C$, where for $n\in \C$, $W_{(n)}$
are the eigenspaces of $L(0)$ with eigenvalues $n$.
A generalized $V$-module 
$W$ is of {\it length $l$}
if there exist 
generalized $V$-submodules $W=W_{1}\supset \cdots 
\supset W_{l+1}=0$ such that $W_{i}/W_{i+1}$ for $i=1, \dots, l$ 
are irreducible $V$-modules. 
A {\it finite length generalized $V$-module} is a generalized $V$-module
of length  $l$ for some $l\in \Z_{+}$. 
{\it Homomorphisms} and {\it isomorphisms} 
between grading-restricted or finite length generalized $V$-modules are 
homomorphisms and isomorphisms between the underlying 
generalized $V$-modules.}
\end{defn}

\begin{rema}
{\rm If $W$ is an $\R$-graded lower-truncated generalized $V$-module
or if $W$ is lower-truncated and generated by one homogeneous element,
then $W$ has a lowest weight. In particular, $V$ or any irreducible 
lower-truncated generalized $V$-module has a lowest weights.}
\end{rema}

\begin{rema}\label{sum-sub-quot}
The category of finite length generalized $V$-modules
is clearly closed under the operation of direct sum, taking 
generalized $V$-submodules and quotient generalized $V$-submodules.
\end{rema}

\begin{prop}\label{contrag}
The contragredient of a generalized $V$-module of length $l$ 
is also of length $l$.
\end{prop}
\pf
Let $W=W_{1}\supset \cdots 
\supset W_{l+1}=0$ be a finite composition series of $W$. 
Then $(W/W_{i})'$ can be naturally embedded into 
$(W/W_{i+1})'$. We view $(W/W_{i})'$ as a generalized $V$-submodule
of $(W/W_{i+1})'$. Then $(W/W_{l+1})'\supset (W/W_{l})'\supset \cdots 
\supset (W/W_{1})'=0$. Moreover $(W/W_{i+1})'/(W/W_{i})'$
is equivalent to $(W_{i}/W_{i+1})'$. Since $W_{i}/W_{i+1}$ is irreducible,
$(W_{i}/W_{i+1})'$ is  irreducible (see \cite{FHL}) and then 
$(W/W_{i+1})'/(W/W_{i})'$ is irreducible. Thus 
$(W/W_{l+1})'\supset (W/W_{l})'\supset \cdots 
\supset (W/W_{1})'=0$ is a composition series of length $l$.
\epfv

\begin{prop}\label{irre-g-r-g-m=>m}
An irreducible grading-restricted generalized $V$-module is a $V$-module.
\end{prop}
\pf
Let $W$ be an irreducible generalized $V$-module. 
Let $W_{(n)}$ be 
the subspace of $W_{[n]}$ containing all eigenvectors 
of $L(0)$ with eigenvalue $n$. Then $\coprod_{n\in \C}W_{(n)}$
is not $0$ and is clearly a $V$-submodule of $W$. 
Since $W$ is irreducible, 
$W=\coprod_{n\in \C}W_{(n)}$, that is,
$W$ is graded by eigenvalues of $L(0)$. 
\epfv

\begin{prop}\label{generators}
A generalized $V$-module of length $l$ 
is generated  by $l$ homogeneous elements whose weights are 
the lowest weights of irreducible $V$-modules. 
\end{prop}
\pf
Let $W$ be a generalized $V$-module of length $l$.
there exist generalized $V$-modules 
$W_{1}\supset \cdots \supset W_{l+1}=0$ such that 
$W_{i}/W_{i+1}$ for $i=1, \dots, l$ 
are irreducible  generalized $V$-modules. 
Since $W_{i}/W_{i+1}$ for $i=1, \dots, l$ 
are irreducible, they are all generated by any nonzero elements.
Let $w_{i}$ for $i=1, \dots, l$
be homogeneous vectors of $W_{i}$ such that 
$w_{i}+W_{i+1}$  are lowest weight 
vectors of 
$W_{i}/W_{i+1}$ for $i=1, \dots, l$, respectively. Since 
$W_{i}/W_{i+1}$ for $i=1, \dots, l$ are irreducible, 
$w_{i}+W_{i+1}$ for $i=1, \dots, l$ are generators of 
$W_{i}/W_{i+1}$.
We claim that $w_{i}$ for $i=1, \dots, l$ 
form a set of generators of $W$. In fact, let $\tilde{W}$ be the 
generalized $V$-submodule generated by $w_{i}$ for $i=1, \dots, l$.
We need to show that $W=\tilde{W}$. Since $W_{l}=W_{l}/W_{l+1}$ is
generated by $w_{l}$, we see that $W_{l}\subset \tilde{W}$. 
Now assume that $W_{m}\subset \tilde{W}$. Then 
since $W_{m-1}/W_{m}$ is generated by $w_{m-1}+W_{m}$,
every element of $W_{m-1}/W_{m}$ is a linear combination of 
elements of the form 
$u^{1}_{n_{1}}\cdots u^{k}_{n_{k}}(w_{m-1}+W_{m})
=u^{1}_{n_{1}}\cdots u^{k}_{n_{k}}w_{m-1}+W_{m}$.
Thus elements of $W_{m-1}$ are linear combinations of elements of 
the form $u^{1}_{n_{1}}\cdots u^{k}_{n_{k}}w_{m-1}+w$ where 
$w\in W_{m}$. Since $u^{1}_{n_{1}}\cdots u^{k}_{n_{k}}w_{m-1}\in \tilde{W}$
and $w\in W_{m}\subset \tilde{W}$, 
$u^{1}_{n_{1}}\cdots u^{k}_{n_{k}}w_{m-1}+w\in \tilde{W}$.
So $W_{m-1}\subset \tilde{W}$. By the principle of induction,
$W=W_{1}\subset \tilde{W}$.
\epfv

\begin{prop}\label{q-fin-dim-gr-res}
A  quasi-finite-dimensional generalized $V$-module is 
grading restricted. An irreducible $V$-module is 
quasi-finite dimensional.
\end{prop}
\pf
If a generalized $V$-module $W$ is quasi-finite dimensional,
then for any $n\in \C$, 
$$\dim W_{[n]}\le \dim \coprod_{\Re{(m)}\le 
\Re{(n)}}W_{[m]}<\infty.$$
If for any $R\in \R$, there exists $n\in \C$ such that 
$\Re{(n)}\le R$ and $W_{[n]}\ne 0$, then clearly 
$\dim \coprod_{\Re{(m)}\le 
R}W_{[m]}=\infty$ for any $R\in \R$. Contradiction. 
So $W$ must also be lower truncated. 

If $W$ is an irreducible $V$-module, then there exists
$h\in \C$ such that
$W=\coprod_{n\in h+\N}W_{(n)}$. Clearly $W$ is quasi-finite 
dimensional.
\epfv

\begin{prop}\label{fin-len-gr-res=>quasi-fin-1}
Every finite-length generalized $V$-module is quasi-finite 
dimensional.
\end{prop}
\pf
Let  $W=W_{1}\supset \cdots \supset W_{n}
\supset W_{n+1}=0$ be a finite composition series of $W$.
Then for $R\in \R$, $\coprod_{\Re{(m)}\le R}
W_{[m]}$ is linearly isomorphic to 
$$\coprod_{i=0}^{n}\coprod_{\Re{(m)}\le R}(W_{i}/W_{i+1})_{[m]}.$$
Since $W_{i}/W_{i+1}$ for $i=0 \dots, n+1$ are 
$V$-modules, that is, they are $L(0)$-semisimple
and grading restricted, by Proposition \ref{q-fin-dim-gr-res},
they are all 
quasi-finite dimensional. So
$$\dim \coprod_{\Re{(m)}\le R}
(W_{i}/W_{i+1})_{[m]}<\infty.$$
Thus $\dim \coprod_{\Re{(m)}\le R}
W_{[m]}<\infty$. 
\epfv

\renewcommand{\theequation}{\thesection.\arabic{equation}}
\renewcommand{\thethm}{\thesection.\arabic{thm}}
\setcounter{equation}{0}
\setcounter{thm}{0}

\section{Cofiniteness conditions and associative algebras}

\begin{defn}\label{cofiniteness-1}
{\rm For a positive integer $n\ge 1$ and 
a weak $V$-module $W$, let $C_{n}(W)$ be the 
subspace of $W$ spanned by elements of the form
$u_{n}w$ where $u\in V_{+}=\coprod_{n\in \Z_{+}}V_{(n)}$ and $w\in W$.
We say that $W$ is {\it $C_{n}$-cofinite} or satisfies the 
{\it $C_{n}$-cofiniteness
condition} if $W/C_{n}(W)$ is finite dimensional. 
In particular, when $n\ge 2$ and $W=V$, we 
say that the vertex operator algebra
$V$ is $C_{n}$-cofinite.}
\end{defn}

\begin{rema}
{\rm The $C_{2}$-cofiniteness condition for the vertex operator algebra
was first introduced
(called ``Condition $C$'') by Zhu
in \cite{Zhu1} and \cite{Zhu2} 
and was used by him to establish the modular invariance 
of the space of characters for the vertex operator algebra.
The $C_{1}$-cofiniteness condition was first introduced 
by Nahm in \cite{N} and was called quasi-rationality there. 
In \cite{L}, generalizing Zhu's $C_{2}$-cofiniteness condition,
Li introduced and studied the $C_{n}$-cofiniteness 
conditions. Note that in the definition above,
when $n=1$ and $W=V$, the cofiniteness condition
is always satisfied. }
\end{rema}

In the case of $n=1$ and $W=V$, 
there is in fact another version of cofiniteness
condition introduced by Li in \cite{L}:

\begin{defn}\label{cofiniteness-2}
{\rm Let $C_{1}^{a}(V)$ be the 
subspace of $V$ spanned by elements of the form 
$u_{n}v$ for $u, v\in V_{+}=\coprod_{n\in \Z_{+}}V_{(n)}$ 
and $L(-1)v$ for $v\in V$. 
The vertex operator algebra $V$ is said to 
be {\it $C^{a}_{1}$-cofinite} or satisfies the 
{\it $C^{a}_{1}$-cofiniteness
condition} if $V/C^{a}_{1}(V)$ is finite dimensional.
For the reason we explain in the remark below, we shall omit the 
superscript $a$ in the notation, that is, we shall 
say $V$ is {\it $C_{1}$-cofinite} or satisfies the 
{\it $C_{1}$-cofiniteness condition} instead of 
$V$ is {\it $C^{a}_{1}$-cofinite} or satisfies the 
{\it $C^{a}_{1}$-cofiniteness
condition}.}
\end{defn}
 
\begin{rema}
{\rm The $C_{1}$-cofiniteness condition in Definition \ref{cofiniteness-2}
can also be defined
for lower-truncated generalized $V$-modules (see \cite{L}). 
But it is now clear that this cofiniteness condition is 
mainly interesting for vertex operator algebras, not for modules,
while the $C_{1}$-cofiniteness condition 
in Definition \ref{cofiniteness-1} is only 
interesting for weak modules, not for vertex operator algebras. 
This is the reason why in the rest of the present paper, we shall omit the 
superscript $a$ (meaning {\it algebra}) in the term ``$C_{1}^{a}$-cofinite''
for $V$, that is, 
when we say that a vertex operator
algebra is $C_{1}$-cofinite, we always mean that it is 
$C_{1}^{a}$-cofinite.}
\end{rema}

\begin{defn}
{\rm A vertex operator algebra $V$ is said
to be of {\it positive energy} if $V_{(n)}=0$ when 
$n<0$ and $V_{(0)}=\C\one$.}
\end{defn}

\begin{rema}
{\rm Positive energy vertex operator algebras are called vertex operator 
algebras of CFT type in some papers, for example, in \cite{GN},
\cite{B} and \cite{ABD}. We use the term ``positive energy''
in this paper
because there are many other
conformal-field-theoretic properties of vertex operator algebras and, 
more importantly, because 
the term ``positive energy'' gives precisely what this condition mean:
If the vertex operator algebra $V$ is the operator product 
algebra of the meromorphic fields of a conformal field theory
so that as an operator acting on this algebra, 
$L(0)=L(0)+\overline{L}(0)$ is equal to the energy
operator, then the energy of 
any state which is not the vacuum is positive and the energy of
the vacuum is of course $0$.}
\end{rema}

\begin{rema}\label{c-n=>c-m}
{\rm Using the $L(-1)$-derivative property, it is 
easy to see that the $C_{n}$-cofiniteness of a 
weak $V$-module $W$ implies the $C_{m}$-cofiniteness 
of $W$ for $1\le m\le n$ and, when $V$ is of positive energy, 
the $C_{2}$-cofiniteness of $V$ implies 
the $C_{1}$-cofiniteness of $V$.}
\end{rema}

The following result is due to Gaberdiel and Neitzke \cite{GN}:

\begin{prop}[\cite{GN}]\label{GN}
Let $V$ be of positive energy and $C_{2}$-cofinite.
Then $V$ is $C_{n}$-cofinite for $n\ge 2$. \epf
\end{prop}


\begin{prop}\label{finite-length-c1}
Assume that all irreducible $V$-modules satisfy the $C_{n}$-cofiniteness 
condition.  Then any finite length generalized module 
is  $C_{n}$-cofinite.
\end{prop}

To prove this result, we need:

\begin{lemma}\label{c1-w2/w1}
For a generalized $V$-module $W_{2}$ and 
a generalized $V$-submodule $W_{1}$ of $W_{2}$, 
$C_{n}(W_{2}/W_{1})=(C_{n}(W_{2})+W_{1})/W_{1}$ as subspaces of 
$W_{2}/W_{1}$.
\end{lemma}
\pf
Note that both $C_{n}(W_{2}/W_{1})$ and 
$(C_{n}(W_{2})+W_{1})/W_{1}$ consists of elements of the form 
$\sum_{i=1}^{k}v^{(i)}_{-n}w^{(i)}+W_{1}$ 
for $v^{(i)}\in V_{+}$
and $w^{(i)}\in W_{2}$. So they are the same.
\epfv

\begin{lemma}\label{w1-c1=>w2-c1}
If $W_{1}$ is a $C_{n}$-cofinite 
generalized $V$-submodule
of a generalized $V$-module $W_{2}$ of finite length such that 
$W_{2}/W_{1}$ is $C_{n}$-cofinite, then 
$W_{2}$ is also $C_{n}$-cofinite. 
\end{lemma}
\pf
Let $X_{1}$ be a subspace of $W_{1}$ such that 
the restrictions to $X_{1}$ of 
the projection from $W_{1}$ to 
$W_{1}/(C_{n}(W_{2})\cap W_{1})$ is a linear isomorphism. 
Let $X_{2}$ be a subspace of $W_{2}$ such that 
the restriction to 
$X_{2}$  of 
the projection from  $W_{2}$ to 
$W_{2}/(C_{n}(W_{2})+W_{1})$ is a linear isomorphism. Then we have 
$W_{2}=C_{n}(W_{2})+X_{1}+X_{2}$ and $W_{2}/C_{n}(W_{2})$ 
is isomorphic to $X_{1}+X_{2}$. 

Since 
$C_{n}(W_{1})\subset (C_{n}(W_{2})\cap W_{1})$, 
$\dim W_{1}/(C_{n}(W_{2})\cap W_{1})\le \dim W_{1}/C_{1}(W_{1})$.
By assumption, $\dim W_{1}/C_{n}(W_{1})<\infty$ and thus 
$\dim W_{1}/(C_{n}(W_{2})\cap W_{1})<\infty$.
By definition $X_{1}$ is 
linearly isomorphic to $W_{1}/(C_{n}(W_{2})\cap W_{1})$.
So $X_{1}$ is also finite-dimensional. 

By definition $X_{2}$ is linearly isomorphic to 
$W_{2}/(C_{n}(W_{2})+W_{1})$ and $W_{2}/(C_{n}(W_{2})+W_{1})$
is linearly isomorphic 
$(W_{2}/W_{1})/((C_{n}(W_{2})+W_{1})/W_{1})$. By Lemma \ref{c1-w2/w1},
$(W_{2}/W_{1})/((C_{n}(W_{2})+W_{1})/W_{1})=
(W_{2}/W_{1})/C_{n}(W_{2}/W_{1})$ which by assumption is 
finite dimensional. Thus
$X_{2}$ is also finite-dimensional. 

Since $W_{2}=C_{n}(W_{2})+X_{1}+X_{2}$ and both $X_{1}$ and $X_{2}$ 
are finite dimensional, $W_{2}/C_{n}(W_{2})$ is finite dimensional. 
\epfv

\noindent {\it Proof of Proposition \ref{finite-length-c1}}\hspace{2ex}
Since $W$ is of finite length, there exist generalized 
$V$-submodules $W= W_{1}\supset \cdots \supset W_{n+1}=
0$ such that $W_{i}/W_{i+1}$ for $i=1, \dots, n$  are 
irreducible. By assumption, $W_{i}/W_{i+1}$ for $i=1, \dots, n$ and
are $C_{n}$-cofinite. Using Lemma \ref{w1-c1=>w2-c1} repeatedly,
we obtain that $W$ is $C_{n}$-cofinite (and in fact $W_{i}$ for
$i=0, \dots, n$ are also $C_{n}$-cofinite). 
\epfv

In the next section, we shall need Zhu's algebra \cite{Zhu1} \cite{Zhu2}
and its generalizations by Dong, Li
and Mason \cite{DLM1} associated to a vertex operator algebra. Here we 
study the relation between the cofiniteness conditions and 
these associative algebras. We first
recall those definitions, constructions and results we need from
\cite{DLM1}. 

For $n\in \N$, define a product $*_{n}$ on $V$ by 
$$u*_{n}v=\sum_{m=0}^{n}(-1)^{m}{m+n\choose n}\res_{x}
x^{-n-m-1}Y((1+x)^{L(0)+n}u, x)v$$
for $u, v\in V$. let $O_{n}(V)$ be the subspace of $V$ 
spanned by elements of the form 
$\res_{x}x^{-2n-2}Y((1+x)^{L(0)+n}u, x)v$
for $u, v\in V$ and of the form
$(L(-1)+L(0))u$ for $u\in V$. 

\begin{thm}[\cite{DLM1}]
The subspace $O_{n}$ is a two-sided ideal of $V$ under the 
product $*_{n}$ and the product $*_{n}$ induces a structure 
of associative algebra on the quotient $A_{n}(V)=V/O_{n}(V)$ with the 
identity $\one +O_{n}(V)$ and with $\omega +O_{n}(V)$ in the center of 
$A_{n}(V)$. \epf
\end{thm}

\begin{rema}
{\rm When $n=0$, $A_{0}(V)$ is the associative algebra first
introduced and 
studied by Zhu in \cite{Zhu1} and \cite{Zhu2}.}
\end{rema}

We shall need the following result 
due to Dong-Li-Mason \cite{DLM2} 
and Miyamoto \cite{Miy}:

\begin{prop}[\cite{DLM2}, \cite{Miy}]\label{c-2n+2=>a-n}
For $n\in \N$, if $V$ is $C_{2n+2}$-cofinite, 
then $A_{n}(V)$ is finite
dimensional. Moreover, $\dim A_{n}(V)\le \dim V/C_{2n+2}(V)$.
\end{prop}
\pf
For the first statement, the case $n=0$ is exactly 
Proposition 3.6 in \cite{DLM2}. 
Here we give a straightforward 
generalization of the proof of Proposition 3.6 in \cite{DLM2}. 
It is slightly different from 
the proof of the general case in the proof of Theorem 2.5 
in \cite{Miy}. Our proof proves the stronger second statement.

By definition, $C_{2n+2}(V)$ are spanned by elements of the 
form $u_{-2n-2}v$ for $u, v\in V_{+}$. Since $V$ is $C_{2n+2}$-cofinite,
there exists a finite dimensional subspace $X$ of $V$ such that 
$X+C_{2n+2}(V)=V$. We need only show that $X+O_{n}(V)=V$. 

By definition, $O_{n}(V)$ is spanned 
by elements of the form 
\begin{eqnarray*}
\lefteqn{\res_{x}x^{-2n-2}Y((1+x)^{L(0)+n}u, x)v}\nn
&&=u_{-2n-2}v+\sum_{k\in \Z_{+}}{\wt u+n\choose k}
u_{k-2n-2}v
\end{eqnarray*}
for $u, v\in V$ and of the form
$(L(-1)+L(0))u$ for $u\in V$. 
We use induction on the weight of elements of $V$.
For any lowest weight vector $w\in V$, 
we have $w=\tilde{w}+\sum_{i=1}^{m}u^{i}_{-2n-2}v^{i}$
where $\tilde{w}\in X$ and $u^{i}, v^{i}\in V$ are homogeneous
for $i=1, \dots, m$. Since $u^{i}_{-2n-2}v^{i}$ are also lowest
weight vectors, $u^{i}_{k-2n-2}v^{i}=0$ for $i=1, \dots, m$ and
$k\in \Z_{+}$.
Then we have 
$$u^{i}_{-2n-2}v^{i}=\res_{x}x^{-2n-2}Y((1+x)^{L(0)+n}u^{i}, x)v^{i}$$
for $i=1, \dots, m$.
We obtain 
$$w=\tilde{w}+\sum_{i=1}^{m}
\res_{x}x^{-2n-2}Y((1+x)^{L(0)+n}u^{i}, x)v^{i}\in X+O_{n}(V).$$
Assume that elements of weights less than $l$ of $V$ are 
contained in $X+O_{n}(V)$. Then for $w$ in $V_{(l)}$, there exists
homogeneous $\tilde{w}\in X$ and homogeneous $u^{i}, v^{i}\in V$ for $i=1, 
\dots, m$ such that $w=\tilde{w}+\sum_{i=1}^{m}u^{i}_{-2n-2}v^{i}$.
Since the weights of $u^{i}_{k-2n-2}v^{i}$ for $i=1, 
\dots, m$ and
$k\in \Z_{+}$ are less than $l$, by induction assumption, 
$$u^{i}_{k-2n-2}v^{i}\in X+O_{n}(V)$$
for $i=1, \dots, m$ and $k\in \Z_{+}$.
Thus 
\begin{eqnarray*}
w&=&\tilde{w}+\sum_{i=1}^{m}u^{i}_{-2n-2}v^{i}\nn
&=&\tilde{w}+\sum_{i=1}^{m}
\res_{x}x^{-2n-2}Y((1+x)^{L(0)+n}u^{i}, x)v^{i}\nn
&& -\sum_{i=1}^{m}\sum_{k\in \Z_{+}}{\wt u^{i}+n\choose k}
u^{i}_{k-2n-2}v^{i}\nn
&\in &X+O_{n}(V).
\end{eqnarray*}
By induction principle, $V=X+O_{n}(V)$, which implies that 
$A_{n}=V/O_{n}(V)$ is finite dimensional.
\epfv

\begin{cor}[\cite{B}]\label{c-2-p-e=>a-n} 
If $V$ is $C_{2}$-cofinite and of positive energy, then for 
$n\in \N$, $A_{n}(V)$ is finite
dimensional.
\end{cor}
\pf
This follow immediately from Proposition \ref{GN} and 
Proposition \ref{c-2n+2=>a-n}.
\epfv

\begin{rema}
{\rm Corollary \ref{c-2-p-e=>a-n} is in fact an easy special case 
of Corollary 5,5 in \cite{B} when the weak $V$-module $M$
there is equal to the vertex operator algebra $V$.}
\end{rema}

\renewcommand{\theequation}{\thesection.\arabic{equation}}
\renewcommand{\thethm}{\thesection.\arabic{thm}}
\setcounter{equation}{0}
\setcounter{thm}{0}

\section{Projective covers of irreducible modules and the finite abelian 
category structure}

\begin{defn}
{\rm Let $\mathcal{C}$ be a full subcategory of generalized $V$-modules. 
A {\it projective object} of $\mathcal{C}$ is an object $W$ of
$\mathcal{C}$ such that for any objects $W_{1}$ and $W_{2}$
of $\mathcal{C}$, any module map $p: W\to W_{2}$ and any surjective 
module map $q: W_{1}\to W_{2}$, there exists a module map 
$\tilde{p}: W\to W_{1}$ such that $q\circ \tilde{p}=p$. 
Let $W$ be an object of $\mathcal{C}$. A {\it projective cover
of $W$ in $\mathcal{C}$} is a projective object $U$ of $\mathcal{C}$
and a surjective module map $p: U\to W$ such that 
for any projective object $W_{1}$ of $\mathcal{C}$ and any 
surjective module map $q: W_{1}\to W$, there exists a surjective 
module map $\tilde{q}: W_{1}\to U$ such that $p\circ \tilde{q}=q$. }
\end{defn}

In general, it is not clear
whether  an object 
of $\mathcal{C}$ has a projective cover 
in $\mathcal{C}$. But we have the following:

\begin{prop} 
If $\mathcal{C}$ is closed under the operations of taking
finite direct sums, quotients and generalized submodules and every
object in $\mathcal{C}$ is completely reducible in $\mathcal{C}$, then
any irreducible generalized $V$-module in $\mathcal{C}$ equipped with
the identity map is a projective cover of the irreducible generalized
$V$-module itself.  
\end{prop} 
\pf 
Let $W$ be an irreducible generalized $V$-module in $\mathcal{C}$ and
$1_{W}: W\to W$ the identity map. We first show that $W$ is projective.
Let $W_{1}$ and $W_{2}$ be objects of $\mathcal{C}$, $p: W\to W_{2}$ a
module map and $q: W_{1}\to W_{2}$ a surjective module map. Since
$W_{2}$ is completely reducible and $W$ is irreducible, $p(W)$ is
irreducible summand of $W_{2}$ and $p$ is an isomorphism from $W$ to
$p(W)$.  Since $W_{1}$ is also completely reducible and $q$ is
surjective, one of the irreducible summand of $W_{1}$ must be isomorphic
to $p(W)$ under $q$.  Let $\tilde{p}: W\to W_{1}$ be the composition of
$p$ and the inverse of the isomorphism from the irreducible summand of
$W_{1}$ above to $p(W)$. By definition, we have $q\circ \tilde{p}=p$. So
$W$ is projective. 

Now let $W_{1}$ be a projective object of 
$\mathcal{C}$ and $q: W_{1}\to W$
a surjective module map. Let $\tilde{q}=q$. 
Then $1_{W}\circ \tilde{q}=q$ and so $(W, 1_{W})$ is 
the projective cover of $W$.
\epfv

In this section, we shall construct projective covers of
irreducible $V$-modules in the category of quasi-finite-dimensional 
generalized $V$-modules when $V$ satisfies certain conditions. 
Our tools are the associative algebras  $A_{n}(V)$, $A_{n}(V)$-modules 
and their relations with generalized $V$-modules. We first need
to recall the constructions and results from \cite{DLM1}.

Let $W$ be a weak $V$-module and let 
$$\Omega_{n}(W)=\{w\in W\;|\;u_{k}w=0\;\mbox{\rm for homogeneous}\;
u\in V, \wt u-k-1\le -n\}.$$

\begin{thm}[\cite{DLM1}]
The map $v\mapsto a_{\swt v-1}$ induces a structure of 
$A_{n}(V)$-module on $\Omega_{n}(W)$. \epf
\end{thm}

The space $\hat{V}$ of operators on $V$ of the form $u_{n}$ for $u\in V$ and 
$n\in Z$, equipped with the Lie Bracket for operators, is a Lie 
algebra by the commutator formula for vertex operators. With the 
grading given by the weights $\wt u-n-1$ of the operators $u_{n}$
when $u$ is homogeneous, $\hat{V}$ is in fact a $\Z$-graded Lie 
algebra. We use $\hat{V}_{(n)}$ to denote the homogeneous subspace of 
weight $n$. Then $\hat{V}_{(0)}$ and $P_{n}(\hat{V})
=\coprod_{k=n+1}^{\infty}\hat{V}_{(-k)}\oplus \hat{V}_{(0)}$ 
are subalgebras of $\hat{V}$. 

\begin{prop}[\cite{DLM1}]
The map given by $v_{\swt v-1} \mapsto v+O_{n}(V)$ is a surjective 
homomorphism of Lie algebras from $\hat{V}_{(0)}$ to
$A_{n}(V)$ equipped with the Lie bracket  induced from 
the associative algebra structure.
\end{prop}

Let $E$ be an $A_{n}(V)$-module. Then it is also a module for 
$A_{n}(V)$ when we view $A_{n}(V)$ as a Lie algebra. By the proposition 
above, $E$ is also a $\hat{V}_{(0)}$-module. Let 
$\hat{V}_{(-k)}$ for $k<n$ act on $E$ trivially. Then 
$E$ becomes a $P_{n}(\hat{V})$-module. Let 
$U(\cdot)$ be the universal enveloping 
algebra functor from the category of Lie algebras to
the category of associative algebras.
Then $U(\hat{V})\otimes_{U(P_{n}(\hat{V}))}E$ 
is a $\hat{V}$-module. If we let elements
of $E\subset U(\hat{V})\otimes_{U(P_{n}(\hat{V}))}E$
to have degree $n$, then the $\Z$-grading on $\hat{V}$ 
induces a $\N$-grading on 
$$U(\hat{V})\otimes_{U(P_{n}(\hat{V}))}E
=\coprod_{m\in N}(U(\hat{V})\otimes_{U(P_{n}(\hat{V})}E)(m)$$
such that $U(\hat{V})\otimes_{U(P_{n}(\hat{V}))}E$ is a graded $\hat{V}$-module. 
By the Poincar\'{e}-Birkhoff-Witt theorem,  
$(U(\hat{V})\otimes_{U(P_{n}(\hat{V}))}E)(m)
=U(\hat{V})_{m-n}E$ for $m\in \Z$,
where $U(\hat{V})_{m-n}$ is the homogeneous 
subspace of $U(\hat{V})$ of degree $m-n$.

For $u\in V$, we define
$Y_{M_{n}}(u, x)=\sum_{k\in \Z}u_{k}x^{-k-1}$. 
These operators give a vertex operator map 
$$Y_{M_{n}}: V\otimes U(\hat{V})\otimes_{U(P_{n}(\hat{V}))}E\to 
U(\hat{V})\otimes_{U(P_{n}(\hat{V}))}E[[x, x^{-1}]$$
for $U(\hat{V})\otimes_{U(P_{n}(\hat{V}))}E$. 
Let $F$ be the subspace of 
$U(\hat{V})\otimes_{U(P_{n}(\hat{V}))}E$ spanned by coefficients of 
$$(x_{2}+x_{0})^{\swt u+n}Y_{M_{n}}(u, x_{2}+x_{0})
Y_{M_{n}}(v, x_{2})w
-(x_{2}+x_{0})^{\swt u+n}Y_{M_{n}}(Y(u, x_{0})v, x_{2})w$$
for $u, v\in V$ and $w\in E$ and 
let 
$$M_{n}(E)=(U(\hat{V})\otimes_{U(P_{n}(\hat{V}))}E)/U(\hat{V})F.$$

\begin{thm}[\cite{DLM1}]\label{dlm-4}
The vector space $M_{n}(E)$ equipped with 
vertex operator map induced from the one for 
$U(\hat{V})\otimes_{U(P_{n}(\hat{V}))}E$
is an $\N$-gradable $V$-module with an $\N$-grading 
$M_{n}(E)=\coprod_{m\in \N}(M_{n}(E))(m)$
induced from the $\N$-grading of $U(\hat{V})\otimes_{U(P_{n}(\hat{V}))}E$ 
such that 
$(M_{n}(E))(0)\ne 0$ and $(M_{n}(E))(n)=E$.
The $\N$-gradable $V$-module satisfies the following universal 
property: For any weak $V$-module $W$ and $A_{n}(V)$-module
map $\phi: E\to \Omega_{n}(W)$, there is a unique 
homomorphism $\bar{\phi}: M_{n}(E)\to W$ of weak $V$-modules
such that $\bar{\phi}((M_{n}(E))(n))=
\phi(E)$.
\end{thm}

\begin{rema}
{\rm In \cite{DLM1}, $M_{n}(E)$ is used to 
denote $U(\hat{V})\otimes_{U(P_{n}(\hat{V}))}E$
and $\bar{M}_{n}(E)$ is used to denote 
what we denote by $M_{n}(E)$ in this paper. 
We use $M_{n}(E)$ instead of $\bar{M}_{n}(E)$
for simplicity. The reader should note
the difference in notations. }
\end{rema}

This finishes our brief discussion of the material in
\cite{DLM1} needed in the present paper. 

The constructions and results
quoted above are for weak modules or $\N$-gradable $V$-modules.
To apply them to our setting, we need the following:

\begin{prop}\label{a-n-1}
If 
$E$ is finite-dimensional, then $M_{n}(E)$
is a generalized $V$-module.
\end{prop}
\pf
By assumption,
$(M_{n}(E))(n)=E$ is finite-dimensional.
Since $L(0)$ preserve the homogeneous subspace $(M_{n}(E))(n)$
of $M_{n}(E)$, we can view $L(0)$
as an operator on the finite-dimensional vector space
$(M_{n}(E))(n)$.
Thus $(M_{n}(E))(n)$ can be decomposed into a direct sum of 
generalized eigenspaces of $L(0)$. This decomposition gives 
$(M_{n}(E))(n)$ a new grading. Since $M_{n}(E)$
is generated by $(M_{n}(E))(n)$, this new grading on 
$(M_{n}(E))(n)$ and the $\Z$-grading on $V$ gives a 
new grading on $M_{n}(E)$ such that the homogeneous 
subspaces are generalized 
eigenspaces of $L(0)$. So $M_{n}(E)$ becomes a
generalized $V$-module. 
\epfv

Recall the $C_{1}$-cofiniteness 
condition for a vertex
operator algebra\footnote{Recall that 
by our convention, the $C_{1}$-cofiniteness condition for a 
vertex operator algebra means the $C_{1}$-cofiniteness condition 
in the sense of \cite{L} or the $C^{a}_{1}$-cofiniteness 
condition.} in the preceding section.

\begin{prop}\label{a-n-1.5}
Let  $V$ be $C_{1}$-cofinite. If a lower-truncated generalized
$V$-module $W$ is finitely generated, 
then $W$ is quasi-finite dimensional. 
\end{prop}
\pf
We need only discuss the case that $W$ is generated 
by one homogeneous element $w$. Since $W$ is lower-truncated,
it is an $\N$-gradable weak $V$-module. 
Since $V$ is $C_{1}$-cofinite,
we know from Theorem 3.10 in \cite{KarL} that 
there are homogeneous $v^{1}, \dots, v^{m}\in V_{+}$ such that 
$W$ is spanned by elements of the form 
$v_{p_{1}}^{i_{1}}\cdots v_{p_{k}}^{i_{k}}w$,
for $1\le i_{1}, \dots, i_{k}\le m$, $p_{1}, \dots, p_{k}\in \Z$ and
$l\in \Z$, satisfying 
$$\wt v^{i_{1}}_{p_{1}}\ge \cdots \ge \wt v^{i_{k_{1}}}_{p_{k_{1}}}>0,$$
$$0>\wt v^{i_{k_{1}+1}}_{p_{k_{1}+1}}\ge \cdots \ge 
\wt v^{i_{k_{2}}}_{p_{k_{2}}}$$
and 
$$\wt v^{i_{k_{2}+1}}_{p_{k_{2}+1}}=\cdots =
\wt v^{i_{k}}_{p_{k}}=0.$$
Since $v^{i_{k_{2}+1}}_{p_{k_{2}+1}}\cdots v^{i_{k}}_{p_{k}}w\in W_{[\swt w]}$,
we see that $W$ is in fact spanned by elements of 
the form $v_{p_{1}}^{i_{1}}\cdots v_{p_{k}}^{i_{k}}\tilde{w}$,
for $1\le i_{1}, \dots, i_{k}\le m$, $p_{1}, \dots, p_{k}\in \Z$,
$l\in \Z$ and $\tilde{w}\in W_{[\swt w]}$, satisfying 
$$\wt v^{i_{1}}_{p_{1}}\ge \cdots \ge \wt v^{i_{k_{1}}}_{p_{k_{1}}}>0$$
and 
$$0>\wt v^{i_{k_{1}+1}}_{p_{k_{1}+1}}\ge \cdots \ge 
\wt v^{i_{k}}_{p_{k}}.$$
But any element $\tilde{w}\in W_{[\swt w]}$ 
satisfies $u^{1}_{j_{1}}\cdots u^{r}_{j_{r}}\tilde{w}=0$ for 
homogeneous $u^{1}, \dots, u^{r}\in V$ and $j_{1}, \dots, 
j_{r}\in \Z$ when $\wt u^{1}_{j_{1}}\cdots u^{r}_{j_{r}}<-n$.
Thus we can take $k=k_{2}$ when $n=0$ and we have
$$\wt v^{i_{1}}_{p_{1}}\ge \cdots \ge \wt v^{i_{k_{1}}}_{p_{k_{1}}}>0$$
and 
$$0>\wt v^{i_{k_{1}+1}}_{p_{k_{1}+1}}\ge \cdots \ge 
\wt v^{i_{k}}_{p_{k}}\ge -n$$
when $n>0$.
{}From these inequalities and the fact that $v^{i_{k_{1}+1}}, 
\dots, v^{i_{k}}\in V_{+}$, we see that when $n>0$,
\begin{eqnarray*}
&k-k_{1}\le n,&\\
&p_{k_{1}+1}, \dots, p_{k}>0,&\\
&p_{k_{1}+1}\le n+\wt v^{i_{k_{1}+1}} -1
\le n+\max(\wt v^{1}, \dots, \wt v^{m})-1&\\
&\cdots,&\\
&p_{k}\le n+\wt v^{i_{k}} -1
\le n+\max(\wt v^{1}, \dots, \wt v^{m})-1.&
\end{eqnarray*}
If 
$\Re{(\wt (v_{p_{1}}^{i_{1}}\cdots v_{p_{k}}^{i_{k}}\tilde{w}))}\le R$,
then we have 
$$\wt v_{p_{1}}^{i_{1}}+\cdots +\wt v_{p_{k_{1}}}^{i_{k_{1}}}
+\wt v^{i_{k_{1}+1}}_{p_{k_{1}+1}}+\cdots +\wt v^{i_{k}}_{p_{k}}
+\Re{(\wt \tilde{w})}\le R.$$
Since 
$$\wt v^{i_{k_{1}+1}}_{p_{k_{1}+1}}+\cdots +\wt v^{i_{k}}_{p_{k}}\ge -n,$$
we obtain
$$0<\wt v_{p_{1}}^{i_{1}}+\cdots +\wt v_{p_{k_{1}}}^{i_{k_{1}}}
+\Re{(\wt \tilde{w})}\le R+n.$$
Combining with the equalities we have above,
we also obtain
$$R+n-\Re{(\wt \tilde{w})}\ge \wt v_{p_{1}}^{i_{1}} \ge \cdots 
\ge \wt v_{p_{k_{1}}}^{i_{k_{1}}}>0.$$
Thus we have 
\begin{eqnarray*}
&k_{1}\le R+n-\Re{(\wt \tilde{w})}=R+n-\Re{(\wt w)},&\\
&p_{1}, \dots, p_{k_{1}}\ge -R-n+\Re{(\wt \tilde{w})}=-R-n+\Re{(\wt w)},&\\
&p_{1}< \wt v^{i_{1}}-1\le \max(\wt v^{1}, \dots, \wt v^{m})-1,&\\
&\cdots,&\\
&p_{k_{1}}< \wt v^{i_{k_{1}}}-1\le \max(\wt v^{1}, \dots, \wt v^{m})-1.&
\end{eqnarray*}
All these inequalities for $k_{1}$, $k-k_{1}$, $p_{1}, \dots, 
p_{k}$ shows that there are only finitely many such numbers and 
thus there are only finitely many elements which span
the subspace $\coprod_{\Re{(l)}\le R}W_{[l]}$ 
of $W$. Thus
the generalized $V$-module $W$ is 
quasi-finite dimensional. 
\epfv

Since the positive energy property and the $C_{2}$-cofiniteness 
condition for $V$
imply the $C_{1}$-cofiniteness condition for $V$,
We have the following consequence:

\begin{cor}
Let $V$ be of positive energy and $C_{2}$-cofinite. Then any 
finitely generated lower-truncated generalized $V$-module
is quasi-finite dimensional. \epfv
\end{cor}

\begin{rema}
{\rm This corollary can also be proved directly using 
a similar argument based on
the spanning set for a weak $V$-module in \cite{B}, without 
using Theorem 3.10 in \cite{KarL}.}
\end{rema}

\begin{cor}\label{a-n-2}
Let  $V$ be $C_{1}$-cofinite. If 
$E$ is finite dimensional, then 
$M_{n}(E)$ is quasi-finite dimensional. 
\end{cor}
\pf
Since $M_{n}(E)$ is a lower-truncated 
generalized $V$-module generated by the finite-dimensional
space $E$, by Proposition \ref{a-n-1.5}, 
it is quasi-finite dimensional. 
\epfv

This result together with Proposition \ref{irre-g-r-g-m=>m}
gives:

\begin{thm}\label{irre-n-g-m=>ord-m}
Let $V$ be $C_{1}$-cofinite.
For an irreducible $\N$-gradable weak $V$-module $W$,
if $\Omega_{0}(W)$ is finite dimensional, then 
$W$ is an ordinary $V$-module. If, in addition, 
$A_{0}(V)$ is semisimple, then 
every irreducible $\N$-gradable 
weak $V$-module is an ordinary $V$-module.
In particular, if $V$ is $C_{1}$-cofinite and 
$A_{0}(V)$ is semisimple, then 
every irreducible lower-truncated generalized $V$-module
is an ordinary $V$-module.
\end{thm}
\pf
Since $\Omega_{0}(W)$ is finite dimensional, 
$M_{0}(\Omega_{0}(W))$ is a quasi-finite dimensional 
generalized $V$-module by Propositions \ref{a-n-1} and \ref{a-n-2}.
The identity map from $\Omega_{0}(W)$ to itself 
extends to a module map from $M_{0}(\Omega_{0}(W))$ 
to $W$. Since $W$ is irreducible, this module map must be 
surjective. Thus $W$ as the image of 
a quasi-finite dimensional 
generalized $V$-module must also be a quasi-finite dimensional 
generalized $V$-module. By Proposition \ref{irre-g-r-g-m=>m},
$W$ must be an ordinary $V$-module. 

If $A_{0}(V)$ is semisimple, 
every irreducible $A_{0}(V)$-module is finite
dimensional. Let $W$ be an irreducible $\N$-gradable weak $V$-module.
Then $\Omega_{0}(W)$ is an irreducible $A_{0}(V)$-module
and hence is finite dimensional. From what we have just proved,
$W$ must be an ordinary $V$-module.
\epfv

The following lemma is very useful:

\begin{lemma}\label{irr-sub-quot}
Let $W$ be a grading-restricted generalized $V$-module
and $W_{1}$ the generalized $V$-submodule of $W$ generated by 
a homogeneous element $w\in \Omega_{0}(W)$. Then there exists 
a generalized $V$-submodule $W_{2}$ of $W_{1}$ such that 
$W_{1}/W_{2}$ is an irreducible $V$-module 
and has the lowest weight $\wt w$.
\end{lemma}
\pf
Since $w\in \Omega_{0}(W)$, $W_{1}$ has the lowest weight $\wt w$. 
Then $(W_{1})_{[\swt w]}$ is an $A_{0}(V)$-module. 
Since $W$ is grading restricted, $(W_{1})_{[\swt w]}$ is 
finite dimensional. 
It is easy to see by induction that any finite-dimensional 
module for an associative algebra is always of finite length. 
In particular, $(W_{1})_{[\swt w]}$ is of finite length. 
Then there exists a $A_{0}(V)$-submodule $M$ 
of $(W_{1})_{[\swt w]}$ such that $(W_{1})_{[\swt w]}/M$
is irreducible. It is also clear that any generalized $V$-submodule 
of $W_{1}$ has a lowest weight. Let $W_{2}$ be the 
sum of all generalized $V$-submodules of $W_{1}$ whose
lowest weight spaces either have weights with real parts 
larger than $\Re{(\wt w)}$ or 
are contained in $M$. Since $W_{2}$ is a sum of generalized $V$-submodule
of $W_{1}$, it is also a generalized $V$-submodule
of $W_{1}$. Since sums of elements of $M$ and elements of
weights with real parts larger than $\Re{(\wt w)}$ cannot be equal to 
$w$, $W_{2}$ cannot contain $w$. Thus $W_{2}$ is a proper 
generalized $V$-submodule of $W_{1}$.

Assume that $W_{3}$ is a proper generalized $V$-submodule of 
$W_{1}$ and contains $W_{2}$ as a 
generalized $V$-submodule. As a generalized submodule of $W_{1}$,
$W_{3}$ also has a lowest weight. Since $M\subset W_{2}\subset W_{3}$,
the lowest weight space of $W_{3}$ must contain $M$. But 
the lowest weight space of $W_{3}$ must be $M$
because otherwise it is an $A_{0}(V)$-module strictly larger than
$M$ but not containing $w$, contradictory to the fact that 
$(W_{1})_{[\swt w]}/M$ is irreducible. We conclude that 
the lowest weight space of $W_{3}$ is in fact $M$. 
Thus, by definition of $W_{2}$, $W_{3}$ is contained in $W_{2}$.
Since by assumption, $W_{3}$ contains $W_{2}$, 
we must have $W_{3}=W_{2}$, proving that $W_{1}/W_{2}$
is irreducible. Since $W$ is grading restricted, 
so is $W_{1}/W_{2}$. Thus $W_{1}/W_{2}$
is an irreducible $V$-module.
Since $w\in W_{1}$ but $w\not \in W_{2}$,
the lowest weight of $W_{1}/W_{2}$ is $\wt w$.
\epfv

Using Proposition \ref{q-fin-dim-gr-res} and 
Lemma \ref{irr-sub-quot}, we obtain the following result:

\begin{prop}\label{q-fin-dim<=>gr-res}
Let $S$ be the set of the lowest weights of irreducible $V$-modules.
Assume that for any $N\in \Z$, 
the set $\{n\in S\;|\;\Re{(n)}\le N\}$
is finite. Then a generalized $V$-module 
is grading-restricted if and only if it is 
quasi-finite dimensional. In particular, the conclusion
holds 
if there are only finitely many inequivalent 
irreducible $V$-modules.
\end{prop}
\pf
In view of Proposition \ref{q-fin-dim-gr-res}, we
need only prove that a grading-restricted generalized
$V$-module is quasi-finite dimensional. 

Let $W$ be 
a  grading-restricted generalized
$V$-module.
We first show that there exists
a subset $S_{0}$ of $S$ such that $W=\coprod_{n\in S_{0}+\N}W_{[n]}$.
Let $w$ be a homogeneous element of $W$. We need only
show that $\wt w\in S+\N$.
The generalized $V$-submodule of $W$ generated by $w$ 
is also grading restricted. Note that this generalized 
$V$-submodule is graded by $\wt w+\Z$. Since it is 
lower truncated and  graded by $\wt w+\Z$, 
it must have a lowest weight of the form $\wt -k$
for some $k\in \N$. 
By Lemma \ref{irr-sub-quot}, this lowest weight $\wt -k$ must be 
the lowest weight of an irreducible $V$-module. 
So we obtain $\wt w-k\in S$ or $\wt w\in S+k\subset S+\N$.

Since for any $N\in \Z$, $\{n\in S_{0}\;|\;\Re{(n)}\le N\}\subset 
\{n\in S\;|\;\Re{(n)}\le N\}$,
by assumption,  $\{n\in S_{0}\;|\;\Re{(n)}\le N\}$
is also finite for any $N\in \Z$. For any $N\in \Z$, 
let $K_{N}$ be a nonnegative integer satisfying
$K_{N}\ge \max \{N-\Re{(n)}\;|\; n\in S_{0}\}$. Then we have 
$$\coprod_{n\in S_{0}+\N, \Re{(n)}\le N}W_{[n]}
\subset \coprod_{i=0}^{K_{N}}\coprod_{n\in S_{0}+i, \Re{(n)}\le N}W_{[n]}.$$
Since for each $n$,
$W_{[n]}$ is finite dimensional and for $i=0, \dots, K_{N}$, 
$\{n\in S_{0}+i\;|\;\Re{(n)}\le N\}$ is finite,
we see that 
$$\coprod_{i=0}^{K_{N}}\coprod_{n\in S_{0}+i, \Re{(n)}\le N}W_{[n]}$$
is finite dimensional. Thus $W$ is quasi-finite dimensional.
\epfv

\begin{prop}\label{quasi-fin=>fin-len}
Assume that there exists $N\in \mathbb{Z}_{+}$ such that 
the real part of the 
lowest weight of any irreducible $V$-module is less than 
or equal to $N$. Then we have:

\begin{enumerate}

\item For any quasi-finite-dimensional generalized 
$V$-module $W$, $\Omega_{0}(W)$ is finite dimensional.

\item Any quasi-finite-dimensional generalized 
$V$-module is of finite length.

\end{enumerate}

\end{prop}
\pf
Let $w$ be a homogeneous element of 
$\Omega_{0}(W)$ and let $W_{1}$  be the generalized 
$V$-submodule of $W$ generated by $w$. 
Then by Lemma \ref{irr-sub-quot}, there exists
a generalized $V$-submodule $W_{2}$ of  $W_{1}$ such that 
$W_{1}/W_{2}$ is irreducible and the lowest weight of 
$W_{1}/W_{2}$ is $\wt w$. By assumption, we have 
$\Re{(\wt w)}\le N$. So $\Omega_{0}(W)\subset \coprod_{\Re{(n)}\le N}
W_{[n]}$. Since $W$ is quasi-finite dimensional, 
$\Omega_{0}(W)$ must be finite dimensional, proving the first 
conclusion.

We now prove the second conclusion. Assume that 
there is a quasi-finite-dimensional generalized $V$-module
$W$ which is not of finite length.
We have just proved that 
$\Omega_{0}(W)$ is finite dimensional. It is easy to see by 
induction
that any finite-dimensional module for an associative algebra
must be of finite length. Since $\Omega_{0}(W)$ is a 
finite-dimensional $A_{0}(V)$-module, there exists 
a finite composition series 
$$M_{0}=\Omega_{0}(W)\supset
M_{1}\supset \cdots\supset M_{k}\supset M_{k+1}=0$$
of $\Omega_{0}$. Take a homogeneous element $w_{1}
\in M_{0}\setminus M_{1}$. 
By Lemma \ref{irr-sub-quot}, we see that $\wt w_{1}$ 
is equal to the lowest weight of an irreducible generalized $V$-module. 

Let $U_{i}$ be the generalized $V$-module generated
by $M_{i}$ for $i=0, \dots, k+1$. Since $W$ is not of finite length,
there exists $i$ such that $U_{i}/U_{i+1}$ is not of finite 
length. Since $U_{i}/U_{i+1}$ is not of finite length, 
it is in particular
not irreducible. So there exists a nonzero proper generalized 
$V$-submodule $U$ of $U_{i}/U_{i+1}$.  
Since $U_{i}/U_{i+1}$ is quasi-finite-dimensional, the 
nonzero proper generalized $V$-submodule $U$ is also 
quasi-finite-dimensional. Then we can repeat the process of obtaining 
the homogeneous element $w_{1}$
above to obtain a nonzero homogeneous element $\tilde{w}\in \Omega_{0}(U)$.

Let $W_{1}=W$ and let $W_{2}$ be the generalized $V$-submodule of $U_{i}$
generated by the elements $w\in U_{i}$ such that $w+W_{i+1}\in 
\Omega_{0}(U)$. Then $W_{2}$ is 
a nonzero proper generalized $V$-submodule of $W_{1}=W$.
We know that $M_{i}/M_{i+1}$ is an irreducible $A_{0}(V)$-submodule
of $\Omega_{0}(U_{i}/U_{i+1})$ and it 
generates $U_{i}/U_{i+1}$. Thus $\Omega_{0}(U)$ cannot
contain any element of $M_{i}/M_{i+1}$. In particular, 
$w_{1}$ cannot be in $W_{2}$.
Let $w_{2}$ be a homogeneous element of $W_{1}$ such that 
$w_{2}+W_{i+1}=\tilde{w}$. 
Since $\tilde{w}\in \Omega_{0}(U)$, by Lemma \ref{irr-sub-quot},
$\wt \tilde{w}$ 
is equal to the lowest weight of an irreducible $V$-module. 
Since $\wt w_{2}=\wt \tilde{w}$, $\wt w_{2}$ is also 
the lowest weight of an irreducible  $V$-module.

Repeating the process above, we obtain an infinite sequence 
$\{W_{i}\}_{i=1}^{\infty}$ of generalized $V$-submodules of $W$ 
such that for $i\in \Z_{+}$, $W_{i+1}$ is a nonzero proper $V$-submodule 
of $W_{i}$ and a sequence 
$\{w_{i}\}_{i=1}^{\infty}$ of homogeneous elements of $W$ such that 
$w_{i}\in W_{i}\setminus W_{i+1}$ and $\wt w_{i}$ is the lowest weights
of an irreducible  $V$-module.

The elements $w_{i}$, $i\in \Z_{+}$, are linearly independent. In fact, if they 
are not. there are $\lambda_{j}\in \C$ and $w_{i_{j}}$ in the 
sequence above for $j=1, \dots, l$ such that $\lambda_{j}$ are not all 
zero, $i_{1}<\cdots<i_{l}$ and
$$\sum_{j=1}^{l}\lambda_{j}w_{i_{j}}=0.$$
We can assume that $\lambda_{1}\ne 0$. Thus $w_{i_{1}}$ 
can be expressed as a linear combination of $w_{i_{2}}, \dots,
w_{i_{l}}$. Since $w_{i}\in W_{i}$ and $W_{i_{j}}\subset W_{i_{2}}$
for $j\ge 2$, we see that $w_{i_{2}}, \dots w_{i_{l}}\in W_{i_{2}}$.
So we see that $w_{i_{1}}$ is a linear combination of elements of 
$W_{i_{2}}$.  So $w_{i_{1}}$ as a linear combination of 
elements of $W_{i_{2}}$ must also be in $W_{i_{2}}$. 
But since $i_{1}>i_{2}$, $W_{i_{2}}$ is a proper submodule 
of $W_{i_{1}}$.
By construction, $w_{i_{1}}\not \in W_{i_{1}+1}\supset W_{i_{2}}$.
Contradiction. So $w_{i}$, $i\in \Z_{+}$, are linearly independent. 

On the other hand, since $\wt w_{i}$ are lowest weights
of irreducible  $V$-modules, the real parts of their weights must be less than 
or equal to $N$. Thus we have a linearly independent infinite subset of 
the finite-dimensional vector space $\coprod_{\Re{(n)}\le N}W_{[n]}$. 
Contradiction. So $W$ is of finite length.
\epfv

\begin{cor}
If there are only finitely many inequivalent 
irreducible $V$-modules, then 
every quasi-finite-dimensional 
generalized $V$-module, or equivalently, every 
grading-restricted generalized $V$-module, is of finite length.
In particular, if $A_{0}(V)$ is finite dimensional, 
every quasi-finite-dimensional 
generalized $V$-module, or equivalently, every 
grading-restricted generalized $V$-module, is of finite length.
\end{cor}
\pf
In this case, a generalized $V$-module is quasi-finite dimensional 
if and only if it is grading restricted by Proposition 
\ref{q-fin-dim<=>gr-res}, and 
the condition in Proposition \ref{quasi-fin=>fin-len}
is clearly satisfied. Thus the conclusion is true.
\epfv

\begin{cor}\label{m-n-e-proj}
Assume that $V$ is $C_{1}$-cofinite and that there exists 
$N\in \Z$ such that the real part of the 
lowest weight of any irreducible $V$-module is less than 
or equal to $N$.
Let $n\in \N$ and let
$E$ be a finite-dimensional $A_{n}(V)$-module.
Then $M_{n}(E)$ is quasi-finite dimensional and  of finite length.
\end{cor}
\pf
Since $V$ is $C_{1}$-cofinite and
$E$ is finite dimensional, by Proposition \ref{a-n-2},
$M_{n}(E)$ is quasi-finite dimensional. By Theorem \ref{quasi-fin=>fin-len},
$M_{n}(E)$ is of finite length.
\epfv

\begin{prop}
Let $E$ be an $A_{n}(V)$-module. Then any eigenspace or generalized 
eigenspace of the operator $\omega+O_{n}(V)\in A_{n}(V)$ on $M$
is also an $A_{n}(V)$-module. 
\end{prop}
\pf
This result follows immediately from the fact that 
$\omega+O_{n}(V)$ is in the center of $A_{n}(V)$.
\epfv

\begin{cor}
Let $E$ be a finite-dimensional $A_{n}(V)$-module,
$\lambda_{1}, \dots, \lambda_{k}$
be the generalized eigenvalues of the operator $\omega+O_{n}(V)$
on $E$ and $E_{(\lambda_{1})}, \dots, E_{(\lambda_{k})}$ the generalized 
eigenspaces of eigenvalues $\lambda_{1}, \dots, \lambda_{k}$,
respectively. Then 
$E_{(\lambda_{i})}$ for $i=1, \dots, k$ are 
$A_{n}$-modules and $E=\coprod_{i=1}^{k}E_{(\lambda_{i})}$. 
\end{cor}

When $M=A_{n}(V)$, we have:

\begin{prop}\label{a_{n}(V)}
Assume that  $A_{n}(V)$ is finite dimensional. Let 
$\lambda_{1}, \dots, \lambda_{k}$
be the generalized eigenvalues of the operator $\omega+O_{n}(V)$
on $A_{n}(V)$. 
Then the generalized eigenspaces 
$(A_{n}(V))_{(\lambda_{i})}$ for $i=1, \dots, k$
are projective $A_{n}(V)$-modules
and $A_{n}(V)=\coprod_{i=1}^{k}(A_{n}(V))_{(\lambda_{i})}$.
\end{prop}
\pf
Since $\omega+O_{n}(V)$ is in the center of $A_{n}(V)$,
$(A_{n}(V))_{(\lambda_{i})}$ is an $A_{n}(V)$-module.
Since $(A_{n}(V))_{(\lambda_{i})}$ is a direct summand of 
the free $A_{n}(V)$-module $A_{n}(V)$ itself, 
it is projective.
\epfv

If an $A_{n}(V)$-module $E$ 
is a generalized eigenspace of the operator $\omega+O_{n}(V)$
with eigenvalue $\lambda$,
we call $E$ a {\it homogeneous $A_{n}(V)$-module of weight 
$\lambda$}. 

\begin{prop}\label{e-weight}
Let $E$ be a homogeneous $A_{n}(V)$-module of weight 
$\lambda$. Then $M_{n}(E)=\coprod_{m\in \N}(M_{n}(E))_{[\lambda -n +m]}$
and 
$(M_{n}(E))_{[\lambda]}=E$,
where $(M_{n}(E))_{[\lambda -n +m]}$ is the generalized eigenspace
of $L(0)$ with eigenvalue $\lambda -n +m$.
\end{prop}
\pf
The operator $L(0)$ acts on $E$ as $\omega+O_{n}(V)$ and 
thus elements of $E$ has weight $\lambda$. The conclusion of
the proposition now follows from the construction of $M_{n}(E)$.
\epfv

\begin{thm}\label{proj-img}
Assume that $V$ is $C_{1}$-cofinite and that there exists 
a positive integer $N$ such that 
$|\Re{(n_{1})}-\Re{(n_{2})}|\le N$ for the lowest
weights $n_{1}$ and $n_{2}$ of any two irreducible $V$-modules.
Let $E$ be a finite-dimensional homogeneous projective $A_{N}(V)$-module
whose weight is equal to the lowest weight of an irreducible $V$-module.
Then $M_{N}(E)$ is projective in the category of 
finite length
generalized $V$-modules.
\end{thm}
\pf
By Propositions \ref{m-n-e-proj}, we know that $M_{N}(E)$
is of finite length. 

Let $W_{1}$ and $W_{2}$ be finite length generalized $V$-modules,
$f: M_{N}(E)\to W_{2}$ a module map and $g: W_{1}\to W_{2}$
a surjective module map. Let the weight of $E$ be $n_{E}$.
Then $n_{E}$ is the lowest weight of some irreducible $V$-module. 
By Proposition \ref{e-weight}, $E$ as a 
subspace of $M_{N}(E)$ is also of weight $n_{E}$. By assumption, 
the set of the real parts of the lowest weights of irreducible $V$-modules 
must be bounded and 
the real parts of weights of any finite length generalized $V$-module 
must be larger than or equal to the real part of the lowest weight of 
an irreducible $V$-module. Thus we see that 
$E$ is in $\Omega_{N}(M_{N}(E))$. 
Then $f(E)$ is in $\Omega_{N}(W_{2})$.
Also $(W_{1})_{(n_{E})}$ and $(W_{2})_{(n_{E})}$
must be in $\Omega_{N}(W_{1})$ and 
$\Omega_{N}(W_{2})$,
respectively. Thus $(W_{1})_{(n_{E})}$ and $(W_{2})_{(n_{E})}$ 
are $A_{N}(V)$-modules. 
Since $f$ and $g$ are module maps
and $g$ is surjective, we have $f(E)\subset (W_{2})_{(n_{E})}$
and $g((W_{1})_{(n_{E})})=(W_{2})_{(n_{E})}$. 
So the restriction $\alpha=f|_{E}: E\to (W_{2})_{(n_{E})}$ and 
$\beta=g|_{(W_{1})_{(n_{E})}}: (W_{1})_{(n_{E})} \to (W_{2})_{(n_{E})}$
of $f$ and $g$ to $E$ and $(W_{1})_{(n_{E})}$, respectively,
are $A_{N}(V)$-module maps and $\beta$
is surjective. Since $E$ is projective, there exists 
an $A_{N}(V)$-module map $\tilde{\alpha}: E\to 
(W_{1})_{(n_{E})}$ such that $\beta\circ \tilde{\alpha}=\alpha$.
By Theorem \ref{dlm-4}, there exists a unique 
$V$-module map $\tilde{f}: M_{N}(E)\to W_{1}$ 
extending $\tilde{\alpha}$.
Since $g\circ \tilde{f}: M_{N}(E)\to W_{2}$ and 
$f: M_{N}(E)\to W_{2}$ are extensions of $\beta\circ \tilde{\alpha}$
and $\alpha$, respectively, and $\beta\circ \tilde{\alpha}=\alpha$,
by the uniqueness in Theorem \ref{dlm-4}, we have 
$g\circ \tilde{f}=f$, proving that $M_{N}(E)$ is 
projective.
\epfv

The following two theorems are the main results of this section:

\begin{thm}\label{proj-img-cor}
Assume that $V$ is $C_{1}$-cofinite and that there exists 
a positive integer $N$ such that 
$|\Re{(n_{1})}-\Re{(n_{2})}|\le N$ for the lowest
weights $n_{1}$ and $n_{2}$ of any two irreducible $V$-modules
and  $A_{N}(V)$ 
is finite dimensional. Then any irreducible 
$V$-module $W$ is has a projective cover in the 
category of  finite length 
generalized $V$-modules.
\end{thm}
\pf
Let $W$ be an irreducible 
$V$-module with the lowest weight $n_{W}$.
Then $W_{(n_{W})}$ is a finite dimensional 
$A_{N}(V)$-module generated by an arbitrary element. 
Since $A_{N}(V)$ 
is finite dimensional, by Proposition \ref{a_{n}(V)},
$A_{N}(V)=\coprod_{i=1}^{k}(A_{N}(V))_{(\lambda_{i})}$
and $(A_{N}(V))_{(\lambda_{i})}$ for $i=1, 
\dots, k$ are projective 
$A_{N}(V)$-modules. Since $W_{(n_{W})}$ 
is an $A_{N}(V)$-module generated by one element,
we have a surjective $A_{N}(V)$-module map
from $A_{N}(V)$ to $W_{(n_{W})}$. But the 
$A_{N}(V)$-module map must preserve the weights,
so under the $A_{N}(V)$-module map, 
$(A_{N}(V))_{(\lambda_{i})}=0$ if $\lambda_{i}\ne 
n_{W}$. Thus we have a surjective $A_{N}(V)$-module map 
from $(A_{N}(V))_{(n_{W})}$ to $W_{(n_{W})}$. 

Now we decompose the finite-dimensional $A_{N}(V)$-module
$(A_{N}(V))_{(n_{W})}$ into a direct sum of 
indecomposable $A_{N}(V)$-modules. Take an
indecomposable $A_{N}(V)$-module $E$ in the decomposition
such that the image of $E$ under the $A_{N}(V)$-module map
from $(A_{N}(V))_{(n_{W})}$ to $W_{(n_{W})}$ 
is not $0$. Since $E$ is a direct summand of the projective 
$A_{N}(V)$-module
$(A_{N}(V))_{(n_{W})}$, $E$ is also projective.
Since $W$ is an irreducible $V$-module, 
$W_{(n_{W})}$ is an irreducible $A_{N}(V)$-module.
Thus the image of $E$ under the $A_{N}(V)$-module map
from $(A_{N}(V))_{(n_{W})}$ to $W_{(n_{W})}$ must be 
equal to $W_{(n_{W})}$. We denote the restriction to $E$ of 
the $A_{N}(V)$-module map
from $(A_{N}(V))_{(n_{W})}$ to $W_{(n_{W})}$ by 
$\alpha$. Then $p$ is surjective.

We first prove that $(E, \alpha)$ is a projective cover of 
$W_{(n_{W})}$ in the category of $A_{N}(V)$-modules.
Let $E_{1}$ be an $A_{N}(V)$-submodule of $E$
such that $E_{1}+\ker \alpha=E$. Let $e_{1}: E_{1}\to E$ be the 
embedding map. Then $\alpha_{1}=p\circ e_{1}$ where 
$\alpha_{1}=\alpha|_{E_{1}}$ is the restriction of $\alpha$ to $E_{1}$. 
Since $\alpha$ is surjective, $\alpha_{1}$ must also be surjective. 
Since $E$ is projective and 
$\alpha_{1}$ is surjective, there exists an $A_{N}(V)$-module
map $\beta_{1}: E\to E_{1}$ such that $\alpha_{1} \circ \beta_{1}=\alpha$. 
Also we have $\alpha_{1} \circ \beta_{1}\circ e_{1}=\alpha\circ e_{1}
=\alpha_{1}$.
Let $E_{2}=\beta_{1}(E_{1})\subset E_{1}$ and $\alpha_{2}=\alpha_{1}|_{E_{2}}
=\alpha|_{E_{2}}:
E_{2}\to W_{(n_{W})}$. Then 
\begin{eqnarray*}
\alpha_{2}(E_{2})&=&(\alpha_{1}\circ \beta_{1})(E_{1})\nn
&=&(\alpha_{1} \circ \beta_{1}\circ e_{1})(E_{1})\nn
&=&\alpha_{1}(E_{1})\nn
&=&W_{(n_{W})},
\end{eqnarray*}
that is, $\alpha_{2}$ is surjective. Since $E$ is projective,
we have an $A_{N}(V)$-module
map $\beta_{2}: E\to E_{2}$ such that $\alpha_{2} \circ \beta_{2}=\alpha$. 
Let $e_{2}: E_{2}\to E$ be the embedding from $E_{2}$ to $E$.
Then we have $p\circ e_{2}=\alpha_{2}$ and so
we also have $\alpha_{2} \circ \beta_{2}\circ e_{2}=\alpha\circ e_{2}=\alpha_{2}$.
Repeating this procedure, we obtain a sequence of 
$A_{N}(V)$-modules $E\supset E_{1}\supset E_{2}\supset
\cdots$ and  $A_{N}(V)$-module maps 
$\beta_{i}: E\to E_{i}$
and $\alpha_{i}: E_{i}\to E$ for $i\in \Z_{+}$ such that 
$E_{i+1}=\beta_{i}(E_{i})$, 
$\alpha_{i} \circ \beta_{i}=\alpha$, $\alpha\circ e_{i}=\alpha_{i}$ and 
$\alpha_{i} \circ \beta_{i}\circ e_{i}=\alpha_{i}$, where 
$e_{i}: E_{i}\to E$ for $i\in \Z_{+}$ are the embedding maps from 
$E_{i}$ to $E$. Since $E$ is of finite length, 
there exists $l\in \Z_{+}$ such that $E_{l+1}=E_{l}$.
Thus $E_{l}=E_{l+1}=\beta_{l}(E_{l})$ and so 
$\beta_{l}\circ e_{l}: E_{l}\to E_{l}$ is surjective. 

We now show that $\gamma=\beta_{l}\circ e_{l}$ must be an isomorphism.
In fact, if not, then $\ker g\ne 0$. Let $K=\ker g$
and $\gamma^{-q}(K)=\gamma^{-1}(g^{-(q-1)}(K))$ for $q\in \N$.
We have a sequence $K\subset \gamma^{-1}(K)\subset \
\gamma^{-2}(K)\subset \cdots$
of $A_{N}(V)$-submodules of $E$. Since 
$E$ is of finite length, there must be $q\in \N$ such that 
$\gamma^{-(q+1)}(K)=\gamma^{-q}(K)$.  Applying $\gamma^{q+1}$ to both sides,
we obtain $K=0$, proving that $\gamma$ is injective. Since $\gamma$ is 
also surjective, 
it is an isomorphism. 

Thus  $(g^{-1}\circ \beta_{l})\circ e_{l}=1_{E_{l}}$, the identity map on 
$E_{l}$. This shows that $E=E_{l}\oplus \ker (g^{-1}\circ \beta_{l})$.
Since $E$ is indecomposable and $E_{l}\ne 0$, we must have 
$\ker (g^{-1}\circ \beta_{l})=0$ and $E=E_{l}$. Thus 
$E_{1}=E$, proving that $(E, \alpha)$ is a projective cover of 
$W_{(n_{W})}$.

By Theorem \ref{proj-img}, $M_{N}(E)$ is a
projective finite length generalized $V$-module and 
by Theorem \ref{dlm-4}, there is a unique module map 
$p: M_{N}(E)\to W$ extending the 
$A_{N}(V)$-module map $\alpha: E\to W_{(n_{W})}$
above. Since $W$ is irreducible and $p\ne 0$, $p$
must be surjective. If there is another projective 
finite length generalized $V$-module $W_{1}$ and a surjective 
module map $q: W_{1}\to W$, there must be a 
module map $\tilde{q}: W_{1}\to M_{N}(E)$
such that $p\circ \tilde{q}=q$. 
Since $(M_{N}(E))_{[n_{W}]}=E$, we see that 
$\tilde{q}((W_{1})_{[n_{W}]})\subset E$. Since 
$q$ is surjective, $q((W_{1})_{[n_{W}]})=W_{(n_{W})}$.
Thus 
\begin{eqnarray*}
\alpha(\tilde{q}((W_{1})_{[n_{W}]}))&=&p(\tilde{q}((W_{1})_{[n_{W}]}))\nn
&=&q((W_{1})_{[n_{W}]})\nn
&=&W_{(n_{W})}.
\end{eqnarray*}
This implies that $\tilde{q}((W_{1})_{[n_{W}]})+\ker \alpha=E$.
Since $(E, \alpha)$ is a projective cover of $W_{(n_{W})}$,
we must have $\tilde{q}((W_{1})_{[n_{W}]})=E$. 
Since $M_{N}(E)$ is generated by $E$,
the image of $W_{1}$ under $\tilde{q}$ is $M_{N}(E)$.
So $\tilde{q}$ is surjective. Thus $(M_{N}(E), p)$
is a projective cover of $W$.
\epfv

Recall from \cite{EO} that an abelian category over $\C$ is called {\it
finite} if every object is of finite length, every space of morphisms is
finite dimensional, there are only finitely many inequivalent simple
objects and every simple object has a projective cover. 

\begin{thm}
Assume that $V$ is $C_{1}$-cofinite and that there exists 
a positive integer $N$ such that 
$|\Re{(n_{1})}-\Re{(n_{2})}|\le N$ for the lowest
weights $n_{1}$ and $n_{2}$ of any two irreducible $V$-modules
and  $A_{N}(V)$ 
is finite dimensional. Then the category of grading-restricted 
generalized $V$-modules is a finite abelian category over $\C$.
\end{thm}
\pf
By Propositions \ref{q-fin-dim<=>gr-res} and \ref{quasi-fin=>fin-len}, 
every object in the category is of finite length.
By Theorem \ref{proj-img-cor}, every object in the category
has a projective cover. Since $A_{N}(V)$ is finite dimensional, there 
are only finitely many inequivalent irreducible (simple) objects. 

Let $W_{1}$ and $W_{2}$ be grading-restricted generalized $V$-modules.
Then they are of finite length. By Proposition \ref{generators},
they are both finitely generated. In particular, $W_{1}$ is generated
by elements of weights whose real parts are less than or equal to
some real number $R$. So module maps from $W_{1}$ to $W_{2}$ 
are determined uniquely by their restrictions to the subspace
$\coprod_{\Re{(n)}\le R}(W_{1})_{[n]}$ of 
$W_{1}$. Since $W_{1}$ and $W_{2}$ are also quasi-finite dimensional
and module maps preserve weights,
these restrictions are maps between finite-dimensional vector spaces
and in particular, the space of these restrictions is finite dimensional. 
Thus the space of module maps from $W_{1}$ to $W_{2}$ is finite 
dimensional.
\epfv

\renewcommand{\theequation}{\thesection.\arabic{equation}}
\renewcommand{\thethm}{\thesection.\arabic{thm}}
\setcounter{equation}{0}

\setcounter{thm}{0}

\section{The tensor product bifunctors and braided tensor
category structure}

We consider a vertex operator algebra 
$V$ satisfying the following conditions:

\begin{enumerate}

\item $V$is $C_{1}$-cofinite\footnote{Recall our convention that 
by a the vertex operator algebra $V$ being $C_{1}$-cofinite, we mean 
that $V$ is $C_{1}$-cofinite in the sense of \cite{L}
or $C^{a}_{1}$-cofinite.}.

\item There exists a positive integer $N$ such that 
$|\Re{(n_{1})}-\Re{(n_{2})}|\le N$ for the lowest
weights $n_{1}$ and $n_{2}$ of any two irreducible $V$-modules
and such that $A_{N}(V)$ is finite dimensional.

\item Every irreducible $V$-module is $\R$-graded and 
$C_{1}$-cofinite\footnote{Recall that by a $V$-module 
being $C_{1}$-cofinite, we mean that the $V$-module 
is $C_{1}$-cofinite in the sense of \cite{H6}
but not necessarily $C^{a}_{1}$-cofinite or $C_{1}$-cofinite in the 
sense of \cite{L}.}.

\end{enumerate}

\begin{prop}\label{c-2+p-e=>cond-1-3}
If $V$ is $C_{2}$-cofinite and of positive energy, 
then $V$ satisfies Conditions 1--3 above.
\end{prop}
\pf
By Remark \ref{c-n=>c-m}, $V$ is $C_{1}$-cofinite.
By Corollary \ref{c-2-p-e=>a-n}, 
$A_{n}(V)$ are 
finite-dimensional for $n\in \N$. In particular, there
are only finitely many inequivalent irreducible $V$-modules.
So Condition 2 is satisfied. 
{}From Theorem 5.10
in \cite{Miy} and Proposition 5.3 in \cite{ABD}, we 
know that every irreducible $V$-modules is 
is $\Q$-graded and $C_{2}$-cofinite
and thus is in particular $C_{1}$-cofinite. 
\epfv

In this section, we shall assume that $V$ satisfies 
Conditions 1--3 above. By Proposition \ref{c-2+p-e=>cond-1-3}, 
if $V$ is $C_{2}$-cofinite and of positive energy, 
this assumption is satisfied.
Thus the results in this section hold if 
$V$ is $C_{2}$-cofinite and of positive energy. 

\begin{prop}
For a vertex operator algebra $V$ satisfying Conditions 1 and 2
above, every irreducible $\N$-gradable 
weak $V$-module is an irreducible $V$-module and 
there are only finitely many inequivalent
irreducible $V$-modules. 
In particular, every lower-truncated irreducible
generalized $V$-module is an irreducible $V$-module. 
\end{prop}
\pf
Since $A_{N}(V)$ is finite dimensional, 
$A_{0}(V)$ as the image of a surjective homomorphism
from $A_{N}(V)$ to $A_{0}(V)$ (see 
Proposition 2.4 in \cite{DLM1})
must also be finite dimensional. Thus there are 
only finitely many inequivalent irreducible $A_{0}(V)$-modules.
By Theorem 2.2.2 in \cite{Zhu2}, there are only finitely
many inequivalent irreducible $V$-modules. 
For an irreducible $\N$-gradable 
weak $V$-module $W$, $\Omega_{0}(W)$ is an irreducible 
$A_{0}(V)$-module by Theorem 2.2.2 in \cite{Zhu2}
and thus must be finite dimensional. 
By Proposition \ref{irre-n-g-m=>ord-m}, every irreducible $\N$-gradable 
weak $V$-module is an irreducible $V$-module. 
\epfv

\begin{prop}
For a vertex operator algebra $V$ satisfying Conditions 1 and 2
above, the category of grading-restricted generalized $V$-modules,
the category of quasi-finite-dimensional 
generalized $V$-modules and the category of finite length  
generalized $V$-modules are the same.
\end{prop}
\pf
{}From Propositions \ref{fin-len-gr-res=>quasi-fin-1}
and \ref{quasi-fin=>fin-len}, we see that 
the category of grading-restricted generalized $V$-modules
is the same as the category of finite length  generalized $V$-modules
and from Propositions \ref{q-fin-dim-gr-res} and 
\ref{q-fin-dim<=>gr-res}, we see that 
the category of grading-restricted generalized $V$-modules
is the same as the category of quasi-finite-dimensional 
generalized $V$-modules.
\epfv

We use $\mathcal{C}$ to denote the category in the proposition above.
In this section, we shall use the results 
obtained in the preceding sections to show
that the category $\mathcal{C}$ satisfies all the assumptions to use the 
logarithmic tensor product theory in \cite{HLZ1} and \cite{HLZ2}. 
Thus by the results of this paper
and the theory developed in \cite{HLZ1} and \cite{HLZ2},
we shall see that $\mathcal{C}$ has a natural structure of 
braided tensor category.

\begin{prop}\label{fusion-rules}
Let $W_{1}$, $W_{2}$ and $W_{3}$ be objects of $\mathcal{C}$.
Then the fusion rule $N_{W_{1}W_{2}}^{W_{3}}$ is finite.
\end{prop}
\pf
By the definition of the category $\mathcal{C}$, 
$W_{1}$, $W_{2}$ and $W_{3}$ are of finite length. 
Since every irreducible $V$-modules is $C_{1}$-cofinite,
by Proposition \ref{finite-length-c1}, 
$W_{1}$, $W_{2}$ and $W_{3}$ are also $C_{1}$-cofinite. 
On the other hand, we also 
know that $W_{1}$, $W_{2}$ and $W_{3}$ are quasi-finite
dimensional. 
Now the proof of Theorem 3.1 in \cite{H6} still works when 
$W_{1}$, $W_{2}$ and $W_{3}$ are quasi-finite-dimensional 
generalized $V$-modules. So the fusion rule 
$N_{W_{1}W_{2}}^{W_{3}}$ is finite.
\epfv

Recall the definition of $W_{1}\hboxtr_{P(z)}W_{2}$
for two generalized $V$-modules $W_{1}$ and $W_{2}$ in
$\mathcal{C}$ in
\cite{HLZ1} and \cite{HLZ2}. Note that
$W_{1}\hboxtr_{P(z)}W_{2}$ depends on our choice 
of $\mathcal{C}$.

\begin{thm}
Let $W_{1}$ and $W_{2}$ be objects in $\mathcal{C}$.
Then $W_{1}\hboxtr_{P(z)}W_{2}$ (defined using the 
category $\mathcal{C}$)
is also in $\mathcal{C}$.
\end{thm}
\pf
By the definition of $W_{1}\hboxtr_{P(z)}W_{2}$, 
it is 
a sum of finite length generalized $V$-submodules 
of $(W_{1}\otimes W_{2})^{*}$. We denote the set of 
these finite length generalized $V$-modules 
appearing in the sum and their finite sums (still
finite length generalized $V$-modules)
by $S$. We want to prove that $W_{1}\hboxtr_{P(z)}W_{2}$
is also of finite length.

Assume that $W_{1}\hboxtr_{P(z)}W_{2}$
is not of finite length. Then take any finite length 
generalized $V$-module $M_{1}$ in $S$. Since $W_{1}\hboxtr_{P(z)}W_{2}$
is not of finite length, $M_{1}$ is not equal to $W_{1}\hboxtr_{P(z)}W_{2}$.
So we can find $M_{2}$ in $S$ such that $M_{1}\subset M_{2}$ but 
$M_{1}\ne M_{2}$. For example, we can take any finite length
generalized $V$-module
in $S$ which is not a generalized submodule of $M_{1}$ and 
then take $M_{2}$ to be the sum of $M_{1}$ and this finite length
generalized $V$-module
in $S$. Since $W_{1}\hboxtr_{P(z)}W_{2}$
is not of finite length, this procedure can continue infinitely
and we get a sequence $\{M_{i}\}_{i\in \Z_{+}}$ 
of finite length generalized $V$-modules in $S$ such that 
$M_{i}\subset M_{j}$ when $i\le j$ but $M_{i}\ne M_{j}$ when 
$i\ne j$. For every $i\in \Z_{+}$, 
since $M_{i}$ is of finite length, $M_{i}/M_{i+1}$ 
is also of finite length. Thus there exists
a generalized $V$-submodule $N_{i}$ of $M_{i}$ such that 
$M_{i+1}\subset N_{i}$ and $(M_{i}/M_{i+1})/(N_{i}/M_{i+1})$ is 
an irreducible $V$-module, or equivalently, 
$M_{i}/N_{i}$ is an irreducible $V$-module. (Note that
by our assumption,
every lower-truncated irreducible generalized $V$-module is an
irreducible $V$-module.) 
Since there are only finitely many equivalence classes
of irreducible $V$-modules, infinitely many 
of the irreducible $V$-modules 
$M_{i}/N_{i}$ for $i\in \Z_{+}$ 
are isomorphic. Let $\{M_{i}/N_{i}\}_{i\in B}$
be an infinite set such that $M_{i}/N_{i}$ for $i\in B\subset \Z_{+}$ are 
isomorphic to an irreducible $V$-module $M$.

By 
Theorem \ref{proj-img-cor},
there exists a projective cover $(P, p)$ of $M$ in the category 
$\mathcal{C}$.
Since $M_{i}/N_{i}$ are isomorphic to 
$M$, there exists surjective module map $\pi_{i}: M_{i}\to M$ whose kernel
is $N_{i}$. Since $P$ is projective,
there exist module maps $p_{i}: P\to M_{i}$ such that 
$\pi_{i}\circ p_{i}=p$. If $p_{i}(P)\subset N_{i}$, then 
$p(P)=\pi_{i}(p_{i}(P))=0$ since the kernel of $\pi_{i}$ is $N_{i}$.
This is contradictory to the surjectivity of $p$. Thus 
$p_{i}(P)$ is not a generalized $V$-submodule of $N_{i}$.

The embedding $M_{i}\to W_{1}\hboxtr_{P(z)}W_{2}$ give
$P(z)$-intertwining maps $J_{i}$ of types ${M_{i}'\choose W_{1}W_{2}}$
as follows: For $m_{i}\in M_{i}$, $w_{1}\in W_{1}$ and $w_{2}\in W_{2}$,
$$\langle m_{i}, J_{i}(w_{1}\otimes w_{2})\rangle
=m_{i}(w_{1}\otimes w_{2}).$$
Let $I_{i}=p_{i}'\circ J_{i}$. Then $I_{i}$ are $P(z)$-intertwining maps
of types ${P'\choose W_{1}W_{2}}$. We now show that these 
intertwining maps are linearly independent. 

Assume that there exist $\lambda_{i}\in \C$, of which only finitely many
are possibly not $0$, such that 
$$\sum_{i\in B}\lambda_{i}I_{i}=0.$$
For $w\in P$, $w_{1}\in W_{1}$ and $w_{2}\in W_{2}$, 
we obtain
\begin{eqnarray*}
0&=&\left\langle w, 
\left(\sum_{i\in B}\lambda_{i}I_{i}\right)(w_{1}\otimes w_{2})\right\rangle\nn
&=&\sum_{i\in B}\lambda_{i}\langle w, 
(p_{i}'(J_{i}(w_{1}\otimes w_{2}))\rangle\nn
&=&\sum_{i\in B}\lambda_{i}\langle p_{i}(w), 
J_{i}(w_{1}\otimes w_{2})\rangle\nn
&=&\sum_{i\in B}\lambda_{i}(p_{i}(w))(w_{1}\otimes w_{2}).
\end{eqnarray*}
Since $w_{1}$ and $w_{2}$ are arbitrary, 
$$\sum_{i\in B}\lambda_{i}p_{i}(w)=0.$$
Since $w\in P$ is also arbitrary, 
$$\sum_{i\in B}\lambda_{i}p_{i}=0.$$
If there exist $i\in B$ such that $\lambda_{i}\ne 0$.
Then there exists $i_{0}\in B$ which is the smallest 
in $B$ such that $\lambda_{i_{0}}\ne 0$. We see that 
$p_{i_{0}}$ can be written as a linear combination of $p_{i}$,
$i\in B$ and $i>i_{0}$. We know that $p_{i}(P)$ is in $M_{i}$
and $p_{i_{0}}(P)$ is in $M_{i_{0}}$ but not in 
$N_{i_{0}}$ which contains but is not equal to $M_{i}$ 
for $i>i_{0}$. Contradiction. So $\lambda_{i}=0$ for all $i\in B$, 
proving the linear independence of $I_{i}$. 

Since $I_{i}$ for $i\in B$ are linearly independent,
the dimension of intertwining maps of type 
${P'\choose W_{1}W_{2}}$ is infinite and thus the fusion rule 
$N_{W_{1}W_{2}}^{P}=\infty$. Since
$W_{1}$, $W_{2}$ and $P$ are in $\mathcal{C}$ and 
so must be quasi-finite dimensional, 
by Proposition \ref{fusion-rules},
$N_{W_{1}W_{2}}^{P}<\infty$. Contradiction.
Thus $W_{1}\hboxtr_{P(z)}W_{2}$ must be of finite length
and thus is in the category $\mathcal{C}$.
\epfv

\begin{cor}
The category $\mathcal{C}$ is closed under $P(z)$-tensor products.\epf
\end{cor}

We now verify the other assumptions in \cite{HLZ1} and \cite{HLZ2}.

\begin{prop}
For any object in $\mathcal{C}$,  the weights 
form a discrete set of rational numbers and there 
exists $K\in \Z_{+}$ such that $(L(0)-L(0)_{s})^{K}=0$
on $W$.
\end{prop}
\pf
By Corollary 5.10 in \cite{Miy}, we know that 
weights of irreducible $V$-modules must be 
in $\Q$. Thus for each irreducible $V$-module,
there exists $h\in \Q$ such that the weights of 
the  irreducible $V$-module are given by $h+k$
for $k\in \N$. For any finite length generalized 
$V$-module $W$, there is a finite composition 
series $W=W_{1}\supset \cdots \supset W_{n}
\supset W_{n+1}=0$
of generalized $V$-submodules of $W$ such that $W_{i}/W_{i+1}$
for $i=1, \dots, n$ are irreducible. Then $W$ as a graded
vector space is 
isomorphic to 
$\coprod_{i=1}^{n}(W_{i}/W_{i+1})$. Since 
$W_{i}/W_{i+1}$ are irreducible, there are $h_{i}\in \Q$
such that the weights of $W_{i}/W_{i+1}$ are
$h_{i}+k$ for $k\in \N$. Thus the weights of 
$W$ are $h_{i}+k$ for $k\in \N$, $i=1, \dots, n$
and clearly form a discrete subset of $\Q$. 

For any weight $m$, $(L(0)-m)W_{[m]}$ is a
subspace of $W_{[m]}$ invariant under $L(0)$. 
It cannot be equal to $W_{[m]}$ since there are 
eigenvectors of $L(0)$ in $W_{[m]}$. Since 
$W/W_{2}$ is an ordinary $V$-module, 
$(L(0)-m)(W/W_{2})_{[m]}=0$, that is, 
$(L(0)-m)W_{[m]}\subset (W_{2})_{[m]}$. 
Similarly we have $(L(0)-m)^{i}(W_{[m]})\subset (W_{i+1})_{[m]}$.
Since $W_{n+1}=0$, we have $(L(0)-m)^{n}(W_{[m]})=0$.
So we can take $K=n$ and then we have 
$(L(0)-L(0)_{s})^{K}=0$ on $W$.
\epfv

Together with Proposition 7.11 in \cite{HLZ2}, we see that 
Assumption 7.10 in \cite{HLZ2} holds for our category $\mathcal{C}$:

\begin{prop}
For any objects $W_{1}$, $W_{2}$ and $W_{3}$ of $\mathcal{C}$,
any logarithmic intertwining operator $\Y$ of type 
${W_{3}\choose W_{1}W_{2}}$ and any $w_{1}\in W_{1}$ and
$w_{2}\in W_{2}$, the powers of $x$ and $\log x$ occurring in 
$\Y(w_{1}, x)w_{2}$ form a unique expansion set 
of the form $D\times \{1, \dots, N\}$ where $D$ 
is a discrete set of rational numbers. \epf
\end{prop}

The definition of $\mathcal{C}$, 
Remark \ref{sum-sub-quot} and Proposition \ref{contrag}
give us the following:

\begin{prop}
The category $\mathcal{C}$ is a full subcategory of 
$\mathcal{G}\mathcal{M}_{sg}$ closed under the operations
of taking contragredients, finite direct sums,
generalized $V$-submodules and quotient generalized $V$-modules.
\epf
\end{prop}

The following result verifies the last assumption
we need:

\begin{thm}
The convergence and the expansion conditions for intertwining maps in
$\mathcal{C}$ hold.  For objects $W_1$, $W_2$, $W_3$, $W_4$, $W_{5}$,
$M_1$ and $M_{2}$ of $\mathcal{C}$, logarithmic intertwining operators
$\mathcal{Y}_{1}$, $\mathcal{Y}_{2}$ and $\mathcal{Y}_{3}$ of types
${W_5}\choose {W_1M_1}$, ${M_1}\choose {W_2M_{2}}$ and ${M_2}\choose
{W_3W_{4}}$, respectively, $z_{1}, z_{2}, z_{3}\in \C$ satisfying
$|z_{1}|>|z_{2}|>|z_{3}|>0$ and $w_{(1)}\in W_{1}$, $w_{(2)}\in
W_{2}$, $w_{(3)}\in W_{3}$, $w_{(4)}\in W_{4}$ and $w'_{(5)}\in
W'_{5}$, the series
$$\sum_{m, n\in {\mathbb C}}\langle w'_{(5)}, \mathcal{Y}_1(w_{(1)}, z_{1})
\pi_{m}(\mathcal{Y}_2(w_{(2)}, z_{2})\pi_{n}(\mathcal{Y}_2(w_{(3)}, z_{3})
w_{(4)}))\rangle_{W_5}$$
is absolutely convergent and can be analytically extended to a 
multivalued analytic function on the region given by 
$z_{1}, z_{2}, z_{3}\ne 0$, $z_{1}\ne z_{2}$, $z_{1}\ne z_{3}$
and $z_{2}\ne z_{3}$ with regular singular points at 
$z_{1}=0$, $z_{2}=0$, $z_{3}=0$, $z_{1}=\infty$, $z_{2}=\infty$, 
$z_{3}=\infty$,
$z_{1}= z_{2}$, $z_{1}= z_{3}$
or $z_{2}= z_{3}$. 
\end{thm}
\pf
By Theorems 11.2 and 11.4 in \cite{HLZ2} (see Remark 12.2
in \cite{HLZ2}), we need only prove that every object of $\mathcal{C}$
satisfies the $C_{1}$-cofiniteness condition, every 
finitely-generated lower-truncated generalized $V$-module
is in $\mathcal{C}$ and every object in $\mathcal{C}$ is quasi-finite
dimensional. Our assumption in the beginning of this section
gives the $C_{1}$-cofiniteness and our definition of the 
category gives the quasi-finite
dimensionality. By Proposition \ref{a-n-1.5}, 
every 
finitely-generated lower-truncated generalized $V$-module
is in $\mathcal{C}$.
\epfv

Using the results above and Theorem 12.13 in \cite{HLZ2}, 
we obtain:

\begin{thm}
The category $\mathcal{C}$, equipped with the tensor
product functor $\boxtimes_{P(1)}=\hboxtr_{P(1)}'$, the unit object $V$,
the braiding, associativity, the left and right unit isomorphisms
given in Subsection 12.2 of \cite{HLZ2}, is a braided tensor 
category. \epf
\end{thm}

\noindent {\small \sc Institut des Hautes \'{E}tudes Scientifiques, 
Le Bois-Marie, 35, Route De Chartres, F-91440 Bures-sur-Yvette, 
France}

\noindent {\it and}

\noindent {\small \sc Department of Mathematics, Rutgers University,
110 Frelinghuysen Rd., Piscataway, NJ 08854-8019 (permanent address)}

\noindent {\em E-mail address}: yzhuang@math.rutgers.edu


\begin{thebibliography}{KWak2}

\bibitem[A]{A} 
T. Abe,
A $\mathbb{Z}\sb 2$-orbifold model of the symplectic fermionic vertex
operator superalgebra, {\em Math. Z.} {\bf 255} (2007), 755--792.

\bibitem[ABD]{ABD}
T. Abe, G. Buhl, C. Dong, Rationality, Regularity,
 and
$C_2$-cofiniteness, {\it Trans. Amer. Math. Soc.} {\bf 356}
(2004), 3391--3402. 

\bibitem[AM1]{AM1} D. Adamovi\'c and A. Milas, Logarithmic 
intertwining operators
and $\mathcal{W}(2,2p-1)$-algebras, {\em Journal of Math. Physics}
{\bf 48}, 073503 (2007).

\bibitem[AM2]{AM2} D. Adamovi\'c and A. Milas, 
On the triplet vertex algebra $\mathcal{W}(p)$, to appear; 
arXiv:0707.1857.

\bibitem[AM3]{AM3} 
D. Adamovi\'c and A. Milas, 
The N=1 triplet vertex operator superalgebras, to appear; 
arXiv:0712.0379.



\bibitem[B]{B}
G. Buhl,  A spanning set for VOA modules, 
{\it J. Algebra} {\bf 254}  (2002),  125--151.

\bibitem[CF]{CF} 
N. Carqueville and M. Flohr, Nonmeromorphic operator
product expansion and $C_2$-cofiniteness for a family of
$\cal{W}$-algebras, {\em J.Phys.} {\bf A39} (2006), 951--966.


\bibitem[DLM1]{DLM1} 
C. Dong, H. Li and G. Mason, Vertex operator algebras and 
associative algebras, {\it J. Algebra} {\bf 206} (1998), 67--96.


\bibitem[DLM2]{DLM2} 
C. Dong, H. Li and G. Mason, Modular
invariance of trace functions in orbifold theory, {\it
Commun. Math. Phys.} {\bf 214} (2000), 1--56.

\bibitem[EO]{EO} 
P. Etingof and V. Ostrik, Finite tensor categories,
{\it Moscow Math. J.} {\bf 4} (2004), 627--654.

\bibitem [FGST1]{FGST1} 
B. L. Feigin, A. M. Ga\u\i nutdinov, 
A. M. Semikhatov, and I. Yu Tipunin, I,
The Kazhdan-Lusztig correspondence for the representation category
of the triplet $W$-algebra in logarithmic conformal field theories
(Russian), {\em Teoret. Mat. Fiz.} {\bf 148} (2006), no. 3, 398--427.

\bibitem [FGST2]{FGST2} 
B. L. Feigin, A. M. Ga\u\i nutdinov, 
A. M. Semikhatov, and I. Yu Tipunin,
Logarithmic extensions of minimal models: characters and modular
transformations, {\em Nucl. Phys.} B {\bf 757} (2006), 303--343.

\bibitem [FGST3]{FGST3} 
B. L. Feigin, A. M. Ga\u\i nutdinov, 
A. M. Semikhatov, and I. Yu Tipunin,
Modular group representations and fusion in logarithmic conformal
field theories and in the quantum group center, {\em Comm. Math.
Phys.} {\bf 265} (2006), 47--93.


\bibitem[FHL]{FHL} I.B. Frenkel, Y.-Z. Huang and J. Lepowsky, On
axiomatic approaches to vertex operator algebras and modules,
preprint, 1989; Memoirs Amer. Math. Soc., Vol. 104, Number 494,
American Math. Soc.  Providence, 1993.

\bibitem[Fl1]{F1}
M. Flohr, On modular invariant partition functions of 
conformal field theories with logarithmic operators, 
{\em Int. J. Mod. Phys.} {\bf A11} (1996), 4147--4172.

\bibitem[Fl2]{F2}
M. Flohr, On fusion rules in logarithmic conformal field theories,
{\em Int. J. Mod. Phys.} {\bf A12} (1996), 1943--1958.

\bibitem[FG]{FG}
M. Flohr and M. R. Gaberdiel, Logarithmic torus 
amplitudes, {\em J. Phys.} {\bf A39} (2006), 1955--1968.

\bibitem[FK]{FK}
M. Flohr and H. Knuth, On Verlinde-Like formulas in $c_{p,1}$
logarithmic conformal field theories, to appear; arXiv:0705.0545.

\bibitem[FGK]{FGK}
M. Flohr, C. Grabow and M. Koehn,
Fermionic Expressions for the characters of $c(p,1)$
logarithmic conformal field theories, {\em Nucl. Phys.} {\bf B768} (2007), 
263--276.

\bibitem[FLM]{FLM} 
I.~B. Frenkel, J.~Lepowsky and A.~Meurman, {\em Vertex
Operator Algebras and the Monster}, Pure and Appl. Math., Vol. 134,
Academic Press, Boston, 1988.

\bibitem[Fu]{Fu} 
J. Fuchs, On nonsemisimple fusion rules and tensor categories, in:
{\it Lie algebras, vertex operator algebras and their applications, 
Proceedings of a conference in honor of James Lepowsky and Robert 
Wilson, 2005}, ed. Y.-Z. Huang and  K. Misra, 
Contemporary Mathematics, Vol. 442, Amer. Math. Soc., Providence, 2007.

\bibitem[FHST]{FHST} J. Fuchs, S. Hwang, A.M. Semikhatov and I. Yu. Tipunin,
Nonsemisimple Fusion Algebras and the Verlinde Formula, {\em Comm.
Math. Phys.} {\bf 247} (2004), no. 3, 713--742.



\bibitem[G]{G} 
V. Gurarie, Logarithmic operators in conformal field
theory, {\em Nucl. Phys.} {\bf B410} (1993), 535--549.

\bibitem[GK1]{GK1}
M. R. Gaberdiel and H. G. Kausch, Indecomposable fusion products,
{\em Nucl. Phys.} {\bf B477} (1996), 298--318.

\bibitem[GK2]{GK2}
M. R. Gaberdiel and H. G. Kausch, A rational logarithmic conformal 
field theory, 
{\em Phys. Lett.} {\bf B386} (1996), 131--137.

\bibitem[GN]{GN}
M. R. Gaberdiel and A. Neitzke, Rationality, quasirationality
and finite W-algebras, {\it Comm. Math.Phys.} {\bf 238} (2003) 305--331.

\bibitem[H1]{H1}
Y.-Z. Huang, A theory of tensor products for module categories for a
vertex operator algebra, IV, {\em J. Pure Appl. Alg.} 100 (1995)
173--216.

\bibitem[H2]{H2} Y.-Z. Huang, Intertwining operator algebras,
genus-zero modular functors and genus-zero conformal field theories,
in: {\em Operads: Proceedings of Renaissance Conferences}, ed. J.-L. Loday,
J. Stasheff, and A. A. Voronov, Contemporary Math., Vol. 202,
Amer. Math. Soc., Providence, 1997, 335--355.

\bibitem[H3]{H3} Y.-Z. Huang, {\em Two-dimensional conformal geometry and
vertex operator algebras}, Progress in Mathematics, Vol. 148,
Birkh\"{a}user, Boston, 1997.

\bibitem[H4]{H4} Y.-Z. Huang, Genus-zero modular functors and
intertwining operator algebras, {\em Internat. J. Math.} 9 (1998), 
845--863.

\bibitem[H5]{H5} Y.-Z. Huang, Riemann surfaces with boundaries and
the theory of vertex operator algebras, in: {\it Vertex Operator
Algebras in Mathematics and Physics}, ed. S. Berman, Y. Billig,
Y.-Z. Huang and J. Lepowsky, Fields Institute Communications, Vol. 39,
Amer. Math. Soc., Providence, 2003, 109--125.


\bibitem[H6]{H6} Y.-Z. Huang, Differential equations and
intertwining operators, {\em Comm. Contemp. Math.} {\bf 7} (2005),
375--400.

\bibitem[H7]{H7} Y.-Z. Huang, Differential equations, duality and
modular invariance, {\it Comm. Contemp. Math.}{\bf 7} (2005), 649--706.

\bibitem[H8]{H8} 
Y.-Z. Huang,  Vertex operator algebras, the Verlinde conjecture 
and modular tensor categories, {\em Proc. Natl. Acad. Sci. USA}
{\bf 102} (2005), 5352--5356. 

\bibitem[H9]{H9} 
Y.-Z. Huang,  Vertex operator algebras, fusion rules and 
modular transformations, in: {\em Non-commutative Geometry and 
Representation Theory in Mathematical Physics}, ed. J. Fuchs, 
J. Mickelsson, G. Rozenblioum and A. Stolin, 
Contemporary Math. Vol. 391, Amer. Math. Soc., Providence, 2005, 
135--148. 

\bibitem[H10]{H10} 
Y.-Z. Huang, Vertex operator algebras and the Verlinde conjecture, 
{\em Comm. Contemp. Math.},  to appear; math.QA/0406291. 

\bibitem[H11]{H11} 
Y.-Z. Huang,  Rigidity and modularity of vertex tensor categories, 
{\em Comm. Contemp. Math.},  to appear; math.QA/0502533. 


\bibitem[HK1]{HK1} 
Y.-Z. Huang and L. Kong, Full field algebras,
{\em Comm. Math. Phys.}{\bf  272}  (2007), 345--396.

\bibitem[HK2]{HK2} 
Y.-Z. Huang and L. Kong, Modular invariance for conformal 
full field algebras, to appear;  arXiv:math/0609570.



\bibitem[HL1]{tensor0}
Y.-Z. Huang and J. Lepowsky, Toward a
theory of tensor products for representations of a vertex operator
algebra, in: {\em Proc. 20th International Conference on Differential
Geometric Methods in Theoretical Physics, New York, 1991},
ed. S. Catto and A. Rocha, World Scientific, Singapore, 1992, Vol. 1,
344--354.

\bibitem[HL2]{tensor-desc}
Y.-Z. Huang and J. Lepowsky, Tensor products of modules for a vertex
operator algebras and vertex tensor categories, in:
     {\em Lie Theory and Geometry,
in honor of Bertram Kostant,}
ed. R. Brylinski, J.-L. Brylinski, V. Guillemin, V. Kac,
Birkh\"{a}user, Boston, 1994, 349--383.

\bibitem[HL3]{tensor1}
Y.-Z. Huang and J. Lepowsky, A theory of tensor products for module
categories for a vertex operator algebra, I, {\em Selecta Mathematica
(New Series)} {\bf 1} (1995), 699--756.

\bibitem[HL4]{tensor2}
Y.-Z. Huang and J. Lepowsky, A theory of tensor products for module
categories for a vertex operator algebra, II, {\em Selecta Mathematica
(New Series)} {\bf 1} (1995), 757--786.

\bibitem[HL5]{tensor3}
Y.-Z. Huang and J. Lepowsky, A theory of tensor
products for module categories for a vertex operator algebra, III,
{\em J. Pure Appl. Alg.} {\bf 100} (1995) 141--171.

\bibitem[HL7]{tensor5}
Y.-Z. Huang and J. Lepowsky, A theory of tensor products for module
categories for a vertex operator algebra, V, to appear.

\bibitem[HLZ1]{HLZ1}
Y.-Z. Huang, J. Lepowsky and L.Zhang, A logarithmic generalization of
tensor product theory for modules for a vertex operator algebra,
{\em Internat. J. Math.} {\bf 17} (2006), 975--1012.

\bibitem[HLLZ]{HLLZ}
Y.-Z. Huang, J. Lepowsky, H. Li and L. Zhang, On the concepts of
intertwining operator and tensor product module in vertex operator
algebra theory, {\em J. Pure Appl. Algebra} {\bf 204} (2006),
507--535.

\bibitem[HLZ2]{HLZ2}
Y.-Z. Huang, J. Lepowsky and L.Zhang, 
Logarithmic tensor product theory for generalized modules for
a conformal vertex algebra, to appear; arXiv:0710.2687.

\bibitem[K1]{K1} H. G. Kausch, Extended conformal algebras
generated by multiplet of primary fields, {\em Phys. Lett.} {\bf 259}
B (1991), 448--455.

\bibitem[K2]{K2} H. G. Kausch, Symplectic fermions, {\em Nucl. Phys.} B {\bf
583} (2000), 513--541.

\bibitem[KarL]{KarL}
M. Karel and H. Li, Certain generating subspaces for vertex operator
algebras, {\em J. Alg.} {\bf 217} (1999),  393--421.

\bibitem[KazL1]{KL1}
D. Kazhdan and G. Lusztig,
Affine Lie algebras and quantum groups,
{\em Duke Math. J., IMRN} {\bf 2} (1991), 21--29.

\bibitem[KazL2]{KL2}
D. Kazhdan and G. Lusztig,
Tensor structures arising {from} affine Lie algebras, I,
{\em J. Amer. Math. Soc.} {\bf 6} (1993), 905--947.

\bibitem[KazL3]{KL3}
D. Kazhdan and G. Lusztig,
Tensor structures arising {from} affine Lie algebras, II,
{\em J. Amer. Math. Soc.} {\bf 6} (1993), 949--1011.

\bibitem[KazL4]{KL4}
D. Kazhdan and G. Lusztig,
Tensor structures arising {from} affine Lie algebras, III, {\em J.
Amer. Math. Soc.} {\bf 7} (1994), 335--381.

\bibitem[KazL5]{KL5}
D. Kazhdan and G. Lusztig,
Tensor structures arising {from} affine Lie algebras, IV,
{\em J. Amer. Math. Soc.} {\bf 7} (1994), 383--453.

\bibitem[Le]{Le} 
J. Lepowsky, From the representation theory of
vertex operator algebras to modular tensor categories in conformal
field theory, commentary on Y.-Z. Huang's PNAS article ``Vertex
operator algebras, the Verlinde conjecture and modular tensor
categories'', {\it Proc. Nat. Acad. Sci. USA} {\bf 102} (2005),
5304--5305.

\bibitem[LL]{LL} 
J. Lepowsky and H. Li, {\em Introduction to Vertex
Operator Algebras and Their Representations}, Progress in Math.,
Birkh\"auser, Boston, 2003.

\bibitem[Li]{L}
H. Li, Some finiteness properties of regular vertex operator algebras,
{\em J. Alg.} {\bf 212} (1999), 495--514.

\bibitem[Mil]{M} A. Milas, Weak modules and logarithmic intertwining
operators for vertex operator algebras, in {\em Recent Developments in
Infinite-Dimensional Lie Algebras and Conformal Field Theory},
ed. S. Berman, P. Fendley, Y.-Z. Huang, K. Misra, and
B. Parshall,
Contemp. Math., Vol. 297,  American Mathematical Society,
Providence,
RI, 2002,  201--225.

\bibitem[Miy]{Miy}
M. Miyamoto, Modular invariance of vertex operator algebras satisfying
$C_{2}$-cofiniteness, {\em Duke Math. J.} {\bf 122} (2004), 51-91.

\bibitem[N]{N}
W. Nahm, Quasi-rational fusion products,
{\em Int. J. Mod. Phys.} {\bf B8} (1994), 3693--3702.

\bibitem[Zha1]{Zha1}
L. Zhang, Vertex operator algebras and Kazhdan-Lusztig's tensor 
category, Ph.D. thesis,
Rutgers University, 2004.

\bibitem[Zha2]{Zha2}
L.Zhang, 
Vertex tensor category structure on a category of
Kazhdan--Lusztig, to appear; math.QA/0701260.

\bibitem[Zhu1]{Zhu1}
Y. Zhu, Vertex operators, elliptic functions and
modular forms, Ph.D. thesis, Yale University, 1990.

\bibitem[Zhu2]{Zhu2}
Y. Zhu, Modular invariance of characters of vertex operator algebras,
{\em J.
Amer. Math. Soc.} {\bf 9} (1996), 237--307.
\end{thebibliography}
\end{document}